\documentclass[12pt]{article}
\usepackage[pdftex]{graphicx}
\usepackage{fullpage}
\usepackage{amsmath}
\usepackage{enumerate}
\usepackage{amssymb}
\usepackage{amsthm}
\usepackage{color} 

\providecommand{\MR}{\relax\ifhmode\unskip\space\fi MR }

\providecommand{\href}[2]{#2}

\begin{document}

\theoremstyle{definition}
\newtheorem*{mydef}{Definition}
\newtheorem*{condition0}{Condition 0}
\theoremstyle{plain}
\newtheorem{theoremsb}{Smooth Billiards Theorem}
\theoremstyle{plain}
\newtheorem*{theoremsi}{Theorem 1:  Smooth Impacts Theorem\label{thm1}}
\title{Smooth Hamiltonian Systems with Soft Impacts}
\author{M. Kloc$^1$, V. Rom-Kedar$^{1,2}$ \\
\normalsize $^1$ Department of Computer Science and Applied Mathematics, \\ \normalsize The Weizmann Institute of Science, Rehovot, Israel \\
\normalsize $^2$ The Estrin Family Chair of Computer Science and Applied Mathematics.
}
\date{\today}
\maketitle

\begin{abstract}
In a Hamiltonian system with impacts (or ``billiard with potential"), a point particle moves about the interior of a bounded domain according to a background potential, and undergoes elastic collisions at the boundaries.  When the background potential is identically zero, this is the hard-wall billiard model.  Previous results on smooth billiard models (where the hard-wall boundary is replaced by a steep smooth billiard-like potential) have clarified how the approximation of a smooth billiard with a hard-wall billiard may be utilized rigorously.  These results are extended  here to models with smooth background potential satisfying some natural conditions. This generalization is then applied to geometric models of collinear triatomic chemical reactions (the models are far from integrable \(n\)-degree of freedom systems with \(n\geq2\)).   The application demonstrates that the simpler analytical calculations for the hard-wall system may be used to obtain qualitative information with regard to the solution structure of the  smooth system and to quantitatively assist in finding solutions of the soft impact system by continuation methods. In particular, stable periodic triatomic configurations are easily located  for the smooth highly-nonlinear two and three degree of freedom geometric models.
\end{abstract}

\section{Introduction}

 The theory of smooth Hamiltonian systems with soft impacts is concerned with a point particle that travels inside a domain  with some given (usually integrable) dynamics, and is repelled from the boundary by a steep smooth potential, essentially as if undergoing an elastic collision with this boundary segment \cite{Dullin98,KozlovBook}.   It is sometimes beneficial to think of the impact system as the limit of smooth soft impact systems.  This approach is particularly useful when studying billiards with inelastic collisions and weak energy dissipation; it is also used as a regularization tool, leading to important persistence results  \cite{KozlovBook}.  The limiting case, in which the soft impacts are replaced by elastic collisions, has been utilized as a model for various physical systems (e.g., a particle moving in a linear gravitational potential or a constant magnetic field  \cite{Berglund00,Berglund96,Dullin98}; chemical reactions \cite{Lerman12}). The dynamics of such systems is non-trivial even in the one-dimensional case \cite{Gorelyshev08}; in higher dimensions, only partial results exist. An extensive theoretical investigation of such systems is presented in \cite{Dullin98}.

Notably, all of the applications mentioned above involve a steep repulsion term, which is replaced in these works by a hard-wall potential for simplicity. Here, we provide conditions under which this approximation is justified, and examples in which it is utilized as a computational tool via continuation methods.
The main results here are Theorems 1 and 2.  In Theorem 1 it is proved,  similarly to the corresponding theorems in \cite{Rapoport07,Turaev98} for the billiard-like potentials, that regular reflections of the smooth impact system are close to those of the hard-wall impact system.  In Theorem 2, we show that the impact Hill's region (which is easily found) may be used to approximate the smooth Hill's region.

Until recently, the soft impact problem was predominantly studied numerically, with no specific attention to the nearly elastic reflections that emerge (for example, classical molecular dynamics simulations involve soft impact problems; see \cite{LeachBook}). In \cite{Lerman12} it was suggested that the techniques developed for steep billiard-like potentials \cite{Rapoport07,Turaev98} could be extended to the soft impact case. Moreover, relying on the current work (i.e., on Theorem 1) and the study of the impact system, qualitative results regarding the behavior of the corresponding soft impact systems were established (such as the existence of homoclinic tangent bifurcations, stable triatomic periodic motion, and, in some cases, simple behavior near the saddle-center point;  see \cite{Lerman12} for details).  Namely, this theorem provides the foundation for analyzing soft impact systems by utilizing the impact limit as demonstrated in \cite{Lerman12} and in section 3 here.  

We note that impact systems belong to the more general field of piecewise smooth dynamical systems, in which the smooth dynamics changes course, possibly undergoing impacts, at some surface.  This field has seen rapid development in the past decade, with particular attention to the new bifurcations that may occur in such systems; see \cite{ChampneysBook, Lamb12} for review, details, and references.

The paper is ordered as follows: in section 2 we first recall a theorem from \cite{Rapoport07,Turaev98} about approximating smooth billiards by hard-wall billiards and then provide its modified formulation for the soft impact case.  This theorem can now be applied to any impact system with a smooth background potential satisfying some natural conditions.   Additionally, we discuss the important appearance of non-trivial  Hill's regions in the new formulation.  In section 3 we present the   application of this theorem to a model of collinear triatomic chemical reactions that was presented in \cite{Lerman12}. We further establish that similar results apply to several (more realistic) extensions of this model and provide numerical simulations of a 3 d.o.f. generalization of this model.  Section 4 provides the proofs of Theorems 1 and 2, and section 5  the summary and discussion.

\section{Formulation of the main result}
We first recall the setup and main relevant results of \cite{Rapoport07}  for the smooth billiard-like dynamics and then modify this set up in order to formulate the new results regarding the smooth impact-like dynamics.\subsection{Billiards and smooth billiard-like potentials}
 Let $D$ be an open bounded subset of $\mathbb{R}^d$ or $\mathbb{T}^d$ with boundary $\partial D=\Gamma_1 \cup \cdots \cup \Gamma_n$, where $\Gamma_i$ are $C^{r+1}$-smooth $(d-1)$-dimensional manifolds of finite area, and $\Gamma^*=\partial \Gamma_1\cup\cdots\cup\partial\Gamma_n$ is the corner set.  The \textit{billiard flow} on $\overline{D}$ is the motion of a point mass in $\overline{D}$ with position $q\in\overline{D}$ and momentum $p\in\mathbb{R}^d$, moving with constant velocity inside $D$ (according to the Hamiltonian $H(q,p)=\frac{p^2}{2}$) and undergoing elastic reflections at $\partial D\setminus \Gamma^*$.  These reflections occur according to the reflection law $p_{out}=p_{in}-2\langle p_{in},n(q)\rangle n(q)$, where $p_{in}$ and $p_{out}$ are the incoming and outgoing velocity vectors, and $n(q)$ is the inward unit normal vector to $\partial D$ at $q$.  The billiard flow can be formally considered as a Hamiltonian system of the form \begin{eqnarray}\label{eq1}
H_b(q,p)=\frac{p^2}{2}+V_b(q), \quad \mbox{where} \quad V_b(q)=\left\{ \begin{array}{cc} 0 & q\in D \\ \mathcal{E} & q\notin D  \end{array} \right.
\end{eqnarray}
for some $\mathcal{E}>0$, and the motion occurs on an energy level $H_b=H^*\in(0,\mathcal{E})$.

Let $V(q;\epsilon)$ be a $C^{r+1}$-smooth \textit{billiard-like potential} : a potential satisfying conditions I-IV listed in Appendix A.   The \textit{smooth billiard flow} is defined by the Hamiltonian \begin{eqnarray}
H(q,p)=\frac{p^2}{2}+V(q;\epsilon),  \quad q(0)\in D, \quad H(q(0),p(0))=H^{*}<\mathcal{E}.
\end{eqnarray}
Roughly, \(V(q;\epsilon)\) is assumed  to tend to zero inside $D$ as $\epsilon\to 0$ and to grow steeply to energies equal to or larger than $\mathcal{E}$ on $\partial D$, thus creating the repulsion from the boundary. Consider a collision point $q_c\in\partial D\setminus \Gamma^*$ belonging to a smooth part of the billiard boundary.  Denote by   \(Q(q;\epsilon)\) a \textit{pattern function}, a function that has a regular limit as \(\epsilon\rightarrow0\) and has the same level sets as the potential \(V(q;\epsilon) \). Choose the coordinates $(x,y)$ so that the hyperplane $x$ is tangent to the potential level surface $Q(q;\epsilon)=Q(q_c;\epsilon)$ (see Conditions I and IIa in Appendix A) and the $y$-axis is the inward normal to this surface at $q=q_c$.

\begin{mydef}
 A reflection at a point $q_c\in\partial D\setminus \Gamma^*$ is called a \textit{regular billiard reflection} if $p_y\neq 0$, and a \textit{non-degenerate tangent billiard reflection} if $p_y=0$ and $p_x^TQ_{xx}p_x>0$.   
\\
 In particular, non-degenerate tangent billiard reflections always satisfy $\|p_x\|\neq 0$.
\end{mydef}

 Under a set of conditions (Conditions I-IV, Appendix A), the billiard flow approximates the smooth billiard flow in the $C^r$-topology for regular reflections and in the $C^0$-topology for non-degenerate tangent reflections (see \cite{Turaev98,Rapoport07} for the two-dimensional and \(n\)-dimensional cases, respectively):
\begin{theoremsb}\cite{Rapoport07}
Let the potential $V(q;\epsilon)$ in the equation \begin{displaymath}
H=\frac{p^2}{2}+V(q;\epsilon)
\end{displaymath} 
satisfy conditions I-IV (Appendix A).  Let $h_{b,t}^\epsilon$ be the smooth billiard flow defined by this equation on an energy surface $H=H^*<\mathcal{E}$, and $b_t$ be the billiard flow in $D$.  Let $\rho_0$ and $\rho_T=b_T\rho_0$ be two inner phase points.  Assume that on the time interval $[0,T]$ the billiard trajectory of $\rho_0$ has a finite number of collisions, and all of them are either regular reflections or non-degenerate tangencies.  Then $\displaystyle h_{b,t}^\epsilon\rho \mathop{\to}_{\epsilon\to 0}b_t\rho$, in the $C^0$ topology for all $\rho$ close to $\rho_0$ and all $t$\ close to $T$.
Finally, if  the billiard trajectory of $\rho_0$ has no tangencies to the boundary on the time interval $[0,T]$,  then $\displaystyle h_{b,t}^\epsilon \mathop{\to}_{\epsilon\to 0}b_t$ in the $C^r$ topology in a small neighborhood of $\rho_0$, and for all $t$ close to $T$.
\end{theoremsb}
\subsection{The soft impact system}
Our main results are  A) extending the above theorem to the soft impact case and B) providing insights regarding the structure of the Hill's region for the soft impact system. 

  Let $U(q)$ be a $C^{r+1}$-smooth potential bounded in the $C^{r+1}$ topology on an open set $\mathcal{D}$ where $\overline{D}\subset \mathcal{D}$, and let $\hat{U}=\min_{q\in\partial D} U(q)$. In particular, since \(|\partial D|\) is finite, $\hat{U}$ is finite.  Let $V(q;\epsilon)$ be a billiard-like potential, as defined in the previous section.  Recall that the minimal barrier height of the billiard-like potential is denoted by \(\mathcal{E}\). To ensure that the particle cannot escape from \(D\), we consider energy levels \(H=H^{*}<\mathcal{E}+\hat U\),and to have motion, we require that  \(\hat U>-\mathcal{E}\):  

\

\noindent\textbf{Condition V} :  $U(q)$ is a $C^{r+1}$-smooth potential bounded in the $C^{r+1}$ topology on an open set $\mathcal{D}$ where $\overline{D}\subset \mathcal{D}$. The minimum $\hat{U}$ of \(U\) on the boundary \(\partial D\) satisfies $\hat U>-\mathcal{E}$.

\

The \textit{impact flow} $\omega_t$ is formally defined as the Hamiltonian flow at an energy level \(H=H^{*}<\mathcal{E}+\hat{U}\) :\begin{eqnarray}
H_{impact}(q,p)=\frac{p^2}{2}+V_b(q)+U(q),\qquad q(0)\in D,
\end{eqnarray}
and similarly, the \textit{smooth impact flow} $h_t^\epsilon$ is the Hamiltonian flow defined by \begin{eqnarray}\label{eq2}
H(q,p)=\frac{p^2}{2}+V(q;\epsilon)+U(q),\qquad q(0)\in D
\end{eqnarray}
on an energy level $H=H^*<\mathcal{E}+\hat{U}$.

\subsubsection{Closeness of impact and smooth trajectories.}

 \
 
We now formulate the extension of the Smooth Billiards Theorem to the smooth impact setting.  In order to do this, the definition of ``non-degenerate tangent reflection" needs to be modified, to account for the curved trajectories under the impact flows:\begin{mydef}
A reflection at a point $q_c\in\partial D\setminus \Gamma^*$ is called a \textit{regular reflection} if $p_y\neq 0$, and a \textit{non-degenerate tangent reflection} if $p_y=0$, $p_x\neq 0$, and $p_x^TQ_{xx}p_x>U_y$.
\end{mydef}

When $U(q)\equiv 0$, this definition coincides with the definition in the billiard systems.  In the billiard case, the non-degeneracy condition excludes nearly-tangent collisions with a concave boundary, whereas here, if $U_y$ is sufficiently negative, such collisions are allowed; see Figure 1.  Notice that regular and non-degenerate tangent collisions must be bounded away from $\partial \mathcal{D}_{Hill}(H^*)\cap \partial\overline D$, where the velocity vanishes.  

\begin{figure}[htbp]
\centering
  \begin{minipage}[b]{5.3 cm}
    \includegraphics[scale=.3]{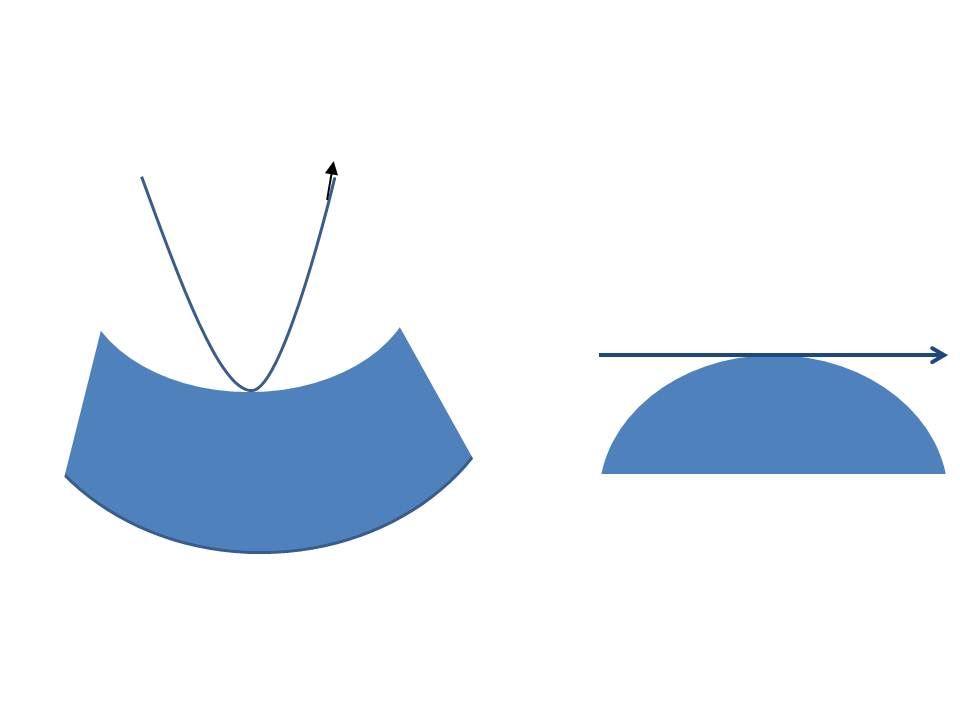}  
  \end{minipage}
\caption{Tangent reflections at concave (left) and convex (right) boundaries.} 
\label{figreflections}
\end{figure}

\begin{theoremsi}
Consider smooth impact systems associated with \begin{displaymath}
H(q,p)=\frac{p^2}{2}+V(q;\epsilon)+U(q), 
\end{displaymath}
and assume that $V(q;\epsilon)$ and $U(q)$ satisfy Conditions I-V (see Appendix A).  Let $\rho_0$ and $\rho_T=\omega_T\rho_0$ be two inner phase points belonging to the energy level $H^*<\mathcal{E}+\hat U$, where $T$ is a finite value.
Provided that on the time interval $[0,T]$ the impact trajectory $\{\omega_t\rho_0:
t\in (0,T]\}$ has a finite number of collisions, and all of them are either regular reflections or non-degenerate tangent reflections, then $\displaystyle h_t^\epsilon\rho\mathop{\to}_{\epsilon\to 0}\omega_t\rho$ in the $C^0$ topology for all $\rho$ close to $\rho_0$ and all $t$ close to $T$.  Moreover, if all collision points are non-tangent, then the above statement is true in the $C^r$ topology.
\end{theoremsi}

The proof is similar to the proof of the soft billiard theorem of \cite{Rapoport07} and is presented in section 4.

\subsubsection{The Hill's region}

The Hill's region of a mechanical Hamiltonian system at a given energy level \(H^{*}\) is defined as the region of allowed motion in the configuration space. This is the region in which the kinetic energy is non-negative, namely, the region in which the potential energy is less than or equal to  \(H^{*}\). The qualitative dependence of the Hill's region on \(H^{*}\)   provides a crude  insight into the highly non-trivial behavior of the flow. Here, we are interested in characterizing

 \begin{displaymath} \mathcal{D^{\epsilon}}_{Hill}(H^*)=
\{q|\;U(q)+V(q;\epsilon)\leq H^*, \quad q\in { \mathcal{D}} \} \end{displaymath} for small \(\epsilon\) by utilizing the impact limit. In the steep billiard-like potential framework (i.e., when \(U(q)\equiv0\)), for any  \(H^{*}\in (0,\mathcal{E})\), the Hill's region of the smooth flow limits,   as \(\epsilon\rightarrow0\), to the billiard domain \(D\). The appearance of the smooth potential $U(q)$ in the impact system  leads to the 
emergence of energy-dependent Hill's regions of the limiting impact system. Let
\begin{equation}
\mathcal{D}^U_{Hill}(H^*)=
\{q|\;U(q)\leq H^*, \quad q\in  \mathcal{D} \} \label{eq:duhill}
\end{equation}   denote the Hill's region of the smooth Hamiltonian \(H=\frac{p^2}{2}+U(q)\), and let us define the Hill's region of the impact system as: \begin{displaymath} \mathcal{D}_{Hill}(H^*)=
\{q|\;U(q)\leq H^*, \quad q\in {\overline D} \}=\mathcal{D}^U_{Hill}(H^*) \cap {\overline D}.\ \end{displaymath} Notice that if the potential attains its global minimum  inside \(D\), and \(H^*\in (\min_{q\in D}U(q),\hat{U})\), the region of allowed motion is strictly inside $D$ (so its boundary is bounded away from \(\partial D\)),
and $\mathcal{D}_{Hill}(H^*)=\mathcal{D}^U_{Hill}(H^*)$. Then, no impacts occur at the energy level $H^*$ and the smooth system trivially limits to the impact system as \(\epsilon\rightarrow0\) by Condition I of Appendix A. Thus, we consider here only the cases when \(H^*\in(\hat{U},\mathcal{E}+\hat{U})\).\footnote{The limit cases where $H^*=\hat{U}$ or $H^*=\hat{U}+\mathcal{E}$ may have interesting behavior. These require further  analysis of collisions with speeds that asymptotically vanish or with escaping orbits. These delicate cases will not be considered here.} 

For this range of energies, the boundary of $\mathcal{D}_{Hill}(H^*)$ is a union of smooth surfaces, some belonging to  \(\partial D\) (denoted by \(\mathcal{\partial D}_{bndry}(H^*)\)) and some belonging to \(\mathcal{D}^U_{Hill}(H^*)\) (denoted by \(\mathcal{\partial D}_{int}(H^*)\)): \begin{equation}
\mathcal{\partial D}_{Hill}(H^*)=\left\{(\partial D\cap\mathcal{D}^U_{Hill}(H^*))\cup(\partial\mathcal{D}^U_{Hill}(H^*)\cap \overline D)\right\}:=\mathcal{\partial D}_{bndry}(H^*)\cup\mathcal{\partial D}_{int}(H^*).
\end{equation}These boundary surfaces meet at the Hill's region corner set \(\Gamma_{Hill}(H^*)=\mathcal{\partial D}_{bndry}(H^*)\cap\mathcal{\partial D}_{int}(H^*)\). Define the  problematic set as the points near the billiard boundary which are also near the Hill's region boundary or in $N(\Gamma^*)$, the set of small neigborhoods around the billiard corners: \(P_{Hill}(H^*)=\left\{ q|\;|U(q)- H^* |<\varsigma \text{ and } |q-{\partial D}|<\xi ֿ\right\}\cup N(\Gamma^*)\) for some small $\varsigma$ and $\xi$; see Figure \ref{proofillustration}. If all of the intersections in the Hill's region corner set are transverse, \(P_{Hill}(H^*)\) is simply an open neighborhood of $\Gamma_{Hill}(H^*)\cup \Gamma^*$, the Hill's region corner set and the billiard corner set.

The main observation here is that away from the problematic set,  for sufficiently small $\epsilon$, the boundary of the impact Hill's region $\mathcal{D}_{Hill}(H^*)$ provides a good smooth approximation to the boundary of the smooth Hill's region $\mathcal{D}^\epsilon_{Hill}(H^*)$. Near transverse intersections of the corner set, the boundaries are \(C^{0}\)-close.  

At non-transverse intersections, we define two types of non-transverse boundary points:  Consider a non-transverse intersection of the corner set at a point $q_c$.  Let $(\bar{x},\bar{y})$ be a fixed local coordinate system with $q_c=(x^*,0)$, so that the hyperplane $\bar{x}$ is tangent to the billiard boundary at \(q_{c}\) and the $\bar{y}$-axis is the inward normal to this boundary at $q_c$. The point \(q_c\) is called \textit{interior  non-transverse boundary point} if $\displaystyle \left.\frac{\partial U}{\partial \bar{y}}\right|_{(x^*,0)}<0$ and a \textit{bifurcating boundary point} if  $\displaystyle \left.\frac{\partial U}{\partial \bar{y}}\right|_{(x^*,0)}\geq0$.

At non-transverse boundary points, the boundaries are also $C^0$-close.  At bifurcating boundary points, one expects more complex behavior, for which bifurcation sequences in $(H,\epsilon)$ space need to be considered.
\\

\noindent\textbf{Theorem 2}:   Assume that the Hill corner region is bounded away from the billiard corner region. Then, for sufficiently small  $\epsilon$:\   
\begin{enumerate}[(i)]
\item Away from the problematic set $P_{Hill}(H^*)$, the boundary of $ \mathcal{D}^\epsilon_{Hill} (H^*)$ is $C^r$-close to the boundary of $\mathcal{D}_{Hill}(H^*)$.
\item Near transverse intersections of the corner set, the smooth Hill's region boundary is   $C^0$-close to the corresponding corner
region of \(\mathcal{\partial D}_{Hill}(H^*)\).
\item    Near interior non-transverse boundary points the smooth Hill's region boundary is   $C^0$-close to  \(\mathcal{\partial D}_{Hill}(H^*)\).
\end{enumerate}

\textbf{Proof}:  See section 4.

\section{Collinear triatomic reaction model}

In this section, we present an example of a smooth impact system for which the computationally-faster hard-wall calculations can be used as a first approximation to the smooth dynamics. This approach allows, for example, to use the impact flow solutions as first guesses for solutions of the smooth system in continuation methods.

  In \cite{Lerman12}, collinear triatomic chemical reactions are modelled as smooth impact systems.   First, we briefly describe the model. Then we present its generalization and detect a new type of stable periodic motion in the generalized system.

For the triatomic collinear reaction $A+BC\to AB+C$, let $r_i$ and $M_i$ $(i\in\{A,B,C\})$  denote the positions and masses, respectively, of the three atoms.  Let $\hat e$ be the unit vector aligned with the molecules.  Since the reaction is collinear, let \\$r_1=(r_A-r_B)\cdot \hat e$ and $ r_2=(r_B-r_C)\cdot \hat e$.  Passing to mass-weighted Jacobi coordinates leads to the Hamiltonian \cite{TannorBook}
\begin{equation}
H(q,p)=\frac{p_1^2}{2}+\frac{p_2^2}{2}+V_r(q_1,q_2),
\label{eq:hreaction}
\end{equation}
where $V_r(q_1,q_2)$ is the potential field, \begin{displaymath}
q_1(r_1,r_2)=\hat a r_1+\hat b r_2\cos \beta, \quad 
q_2(r_2)=\hat b r_2\sin\beta,
\end{displaymath}
\begin{displaymath}
\hat a=\sqrt{\frac{M_A(M_B+M_c)}{M_A+M_B+M_C}}, 
\quad \hat b=\sqrt{\frac{M_c(M_B+M_A)}{M_A+M_B+M_C}},
\end{displaymath}
and \begin{displaymath}
\beta=\arccos \sqrt{\frac{M_AM_C}{(M_A+M_B)(M_B+M_C)}}.
\end{displaymath}

The potential \(V_r(q_1,q_2)\) of eq (\ref{eq:hreaction}) is just the potential energy surface (PES) at the collinear configuration. Much effort has been dedicated in recent years to finding a good form for the PES  \cite{LawleyBook}.

The reaction region is defined to be the region where both \(r_{1,2}\) are bounded, whereas large \(r_{1}\) and bounded \(r_{2}\) (respectively, large  \(r_{2}\) and bounded \(r_{1}\)) correspond to the reactant (respectively, product) channel. A trajectory with initial conditions in the reactant channel enters the reaction region and then may exit  through either channel. If it exits through the product channel, the reaction is realized. Finding the reaction rates analytically, even  in this highly simplified model, is practically impossible as it is typically chaotic. Transition state theory attempts to approximate these rates by examining the local dynamics near a saddle point of the potential, and is known to be problematic. One should also note that relating the numerically-calculated reaction rates of  such a single molecular reaction model (classical, semi-classical or quantum, for collinear or for the three dimensional reaction model) to the kinematic reaction rates is also an open problem \cite{Wu06}. Nonetheless, understanding the qualitative features of the collinear model may help to shed light on these challenging open problems.

Thus, here, as in   \cite{Lerman12}, we concentrate on qualitative features of the dynamics, without addressing the practical aspects of reaction rate calculations.  In   \cite{Lerman12} it was suggested that in triatomic collinear reactions in which there is a single unstable triatomic configuration, the level sets of  \(V_r(q_1,q_2)\) may be modeled by the following potential form that  has three components, each having a distinct geometrical meaning:  \begin{equation}
H(q,p; \epsilon)=\frac{p^2}{2}+bV_b(q)+aV_a(q)+cV_{farfield}(q).\label{eq:hgeometricqp}
\end{equation}
The first term, $bV_b$, corresponds to the strong nuclear repulsion of the diatoms at small distances (at very small \(r_{i},i=1,2\)), and is thus modeled by a smooth billiard-like potential. The second term, $aV_a$, represents the local  interactions near the transition state (the unstable triatomic configuration), and is modeled in    \cite{Lerman12}  by the quadratic potential \begin{eqnarray}
aV_a(q)=\frac{1}{2}(q-q_s)^TA(q-q_s),
\end{eqnarray}
where $A$ is a symmetric 2x2 matrix with negative determinant. Then, the unstable triatomic configuration corresponds to a saddle-center fixed point in the phase space.  The last term,  $cV_{farfield}$, is assumed to be small in the reaction region, and is chosen to have the correct asymptotic form at the reactant and product channels (large \(r_1\) or large \(r_2\)).  Since the main interest here and in     \cite{Lerman12}  is in studying the dynamics in the reaction region, no specific form for this term is needed here.

The Hamiltonian (\ref{eq:hgeometricqp}) is put into normal form by rotating the $(q,p)$ coordinates into  the coordinates that are aligned with the center and the saddle subspaces, so that in these coordinates the saddle-center (the reaction barrier) linear part is diagonal. Then, the smooth Hamiltonian  model is of the form: \begin{displaymath}
H_{smooth}(u_1,u_2,v_1,v_2)=\frac{v_1^2}{2}+\frac{v_2^2}{2}+U(u_1,u_2)+V_b(u_1,u_2;\epsilon)
\end{displaymath}
with quadratic background potential \begin{displaymath}
U(u_1,u_2)=\frac{\omega^2}{2}(u_1-u_{1s})^2-
\frac{\lambda^2}{2}(u_2-u_{2s})^2,
\end{displaymath}
and the billiard-like potential is taken here to be symmetric in the wedge and of the form:\begin{displaymath}
V_b(u_1,u_2; \epsilon)=b\cdot \exp\left(\frac{-(u_1\sin(\beta/2)+u_2\cos(\beta/2))}{\epsilon}\right)
+b\cdot \exp\left(\frac{-(u_1\sin(\beta/2)-u_2\cos(\beta/2))}{\epsilon}\right).
\end{displaymath}
 In \cite{Lerman12} this geometrical model (in a  generalized asymmetric form) was introduced and several results regarding the motion at energies close to the barrier energy (such as the existence of homoclinic bifurcations and, consequently, of elliptic islands near some asymmetric configurations) were established by utilizing the impact limit and relying on the smooth impact theorem (Theorem 1 of section 2).  Figure 2 demonstrates the validity of Theorem 2.  In particular, it is seen that when $\mathcal{D}_{Hill}(H^*)$ has corners, the smooth Hill's region converges towards the singular domain in the $C^0$ topology.

 Here we consider the behavior near a particular type of symmetric periodic orbit (see Figure \ref{figtraj}) which is far from the saddle-center point. The existence and stability of such an orbit is established first in the impact system and then, by continuation,  in the smooth case. In fact, we consider here a higher-dimensional version of the geometric model by adding a  ``group of oscillators" -- a standard extension used in chemistry to reflect the effect of small oscillatory modes \cite{BialekBook}. Specifically, we add \(n-2\) separable nonlinear oscillators and a weak coupling term:
  \begin{displaymath}
H_{smooth}(u_1,...,u_n,v_1,...,v_n)=\frac{v_1^2}{2}+\frac{v_2^2}{2}+U(u_1,u_2)+V_b(u_1,u_2;\epsilon)
+\sum_{k=3}^n \left( \frac{v_k^2}{2}+U_k(u_k) \right)+\delta U_{coup}(u_1,...,u_n).
\end{displaymath}
In the numerical computations we  let $n=3$,   take a quartic symmetric potential \begin{displaymath}
U_3(u_3)=\frac{(\kappa\omega)^2}{2}(u_3-u_{3s})^2+\frac{1}{4}(u_3-u_{3s})^4,
\end{displaymath}
with $\kappa$ chosen to be irrational, and define\begin{displaymath}
U_{coup}(u_1,u_2,u_3)=\sin(u_1-u_2)+\sin(u_2-u_3)+\sin(u_3-u_1).
\end{displaymath}
First let $\delta=0$.  In this case, $u_3$ and $v_3$ do not affect the behavior of the system, so we consider only $u_1$, $u_2$, $v_1$, and $v_2$.
By the smooth impact theorem, we can approximate this system by the impact flow, where the steep smooth billiard-like potential is replaced by a hard-wall billiard in a  wedge with  upper (respectively, lower) boundary defined by the unit vector $[\cos\beta/2, \quad \sin\beta/2]^T$ (respectively,  $[\cos\beta/2, \quad -\sin\beta/2]^T$).  The particle undergoes elastic collisions at the upper billiard boundary according to the reflection law \begin{eqnarray*}
v_1&\mapsto& v_1\cos\beta+v_2\sin\beta\\
v_2&\mapsto& v_1\sin\beta-v_2\cos\beta 
\end{eqnarray*} and at the lower billiard boundary according to \begin{eqnarray*}
v_1&\mapsto& v_1\cos\beta- v_2\sin\beta\\
v_2&\mapsto& -v_1\sin\beta-v_2\cos\beta.
\end{eqnarray*}
In the interior, the integrable linear system \begin{eqnarray*}
\frac{du_1}{dt}=v_1, \quad \frac{du_2}{dt}=v_2, \quad \frac{dv_1}{dt}=-\omega^2(u_1-u_{1s}), \quad \frac{dv_2}{dt}=\lambda^2(u_2-u_{2s})
\end{eqnarray*}
can be solved analytically.
Here we consider the symmetric case, where the saddle point lies on the $u_1$-axis, so \(u_{2s}=0\).      Then the linear flow is \begin{eqnarray*}
 u_1(t)&=& u_{1s}+(u_{10}-u_{1s})\cos(\omega t)+\frac{v_{10}}{\omega}\sin(\omega t)\\
 u_2(t)&=& u_{20}\cosh(\lambda t)+\frac{v_{20}}{\lambda}\sinh(\lambda t)\\
 v_1(t)&=& -\omega(u_{10}-u_{1s})\sin(\omega t)+v_{10}\cos(\omega t)\\
 v_2(t)&=& \lambda u_{20}\sinh(\lambda t)+v_{20}\cosh(\lambda t).
 \end{eqnarray*} 

Consider initial conditions on the \(u_{1}\) axis, to the right of the saddle point, with positive initial velocity only in the perpendicular  \(u_{2}\) direction (i.e., $u_{10}>u_{1s},u_{20}=0,v_{10}=0,v_{20}>0)$.  For every such initial position $(u_{10}, 0)$, we propose that by the symmetry of the problem, there is some $v_{20}$ such that the trajectory with initial conditions $(u_{10}, 0, 0, v_{20})$ is 2-periodic:  it hits the upper boundary at a right angle.\footnote{In fact, using $t_c$ as a parameter, one can find an implicit equation for the $v_{20}$ for which such a trajectory is 2-periodic.  It is not difficult to see that a discrete family of $t_c$ values may be defined in this way; however, we expect that these solutions correspond to a period-two orbit only for values of $u_{10}$ that are exponentially close to $u_{1s}$.}   Below, we find such orbits for $\epsilon=0$, show that these orbits can be continued for $\epsilon>0$, and determine their stability.  An example of such an orbit for $\epsilon=0$ and $\epsilon=0.1$ is shown in Figure \ref{figtraj}.

\begin{figure}[htbp]
\centering
  \begin{minipage}[b]{9 cm}
    \includegraphics[scale=.5]{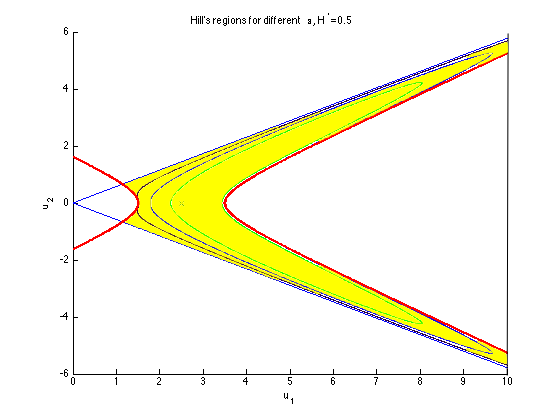}  
  \end{minipage}
\caption{Demonstration of the Hill's regions' convergence (Theorem 2):  The impact Hill's region $\mathcal{D}_{Hill}(H^*)$ is shown in yellow.  The corresponding Hill's region boundaries (where $U(q)+V(q;\epsilon)=H^*$) for $\epsilon=0.1$ (black), $\epsilon=0.2$ (blue), $\epsilon=0.3$  (green), are shown.  For $\epsilon=.001$ and $\epsilon=.01$, the smooth and impact regions are indistinguishable.  The $C^0$ convergence at the corners is apparent.  Parameters: $H^*=0.5$, $\omega=1$, $\lambda=\sqrt{2}$, $\beta=\pi/3$, $u_{1s}=2.5$, $u_{2s}=0$, $b=10$.   } \label{hill}
\end{figure}

\begin{figure}[htbp]
\centering
  \begin{minipage}[b]{5.3 cm}
    \includegraphics[scale=.3]{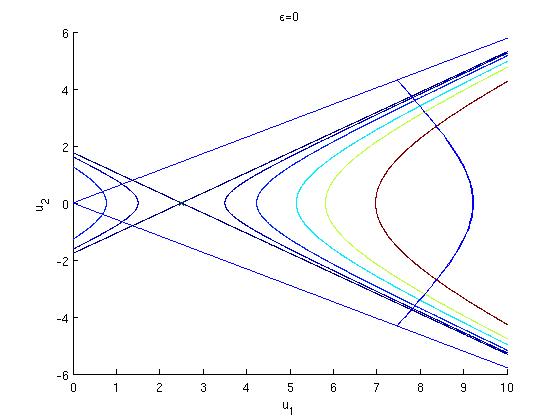}  
  \end{minipage}
  \begin{minipage}[b]{5.3 cm}
    \includegraphics[scale=.3]{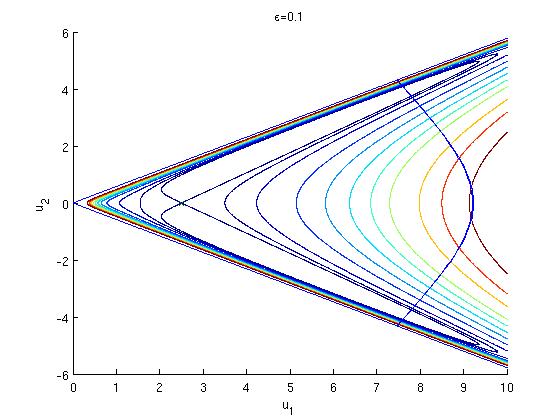}  
  \end{minipage}
\caption{A  2-periodic trajectory for $\epsilon=0$ (with the level lines of the background potential $U(q)$) and $\epsilon=0.1$ (with the level lines of the full potential $U(q)+V(q;\epsilon)$).  Parameters: $\omega=1$, $\lambda=\sqrt{2}$, $\beta=\pi/3$, $u_{1s}=2.5$, $u_{2s}=0$, $b=10$.   The trajectory shown with  $\epsilon=0$ has \(u_{10}=9.23,v_{20}=4.89\) and corresponding energy \(H=35.064\), whereas at $\epsilon=0.1$ we set  \(u_{10}=9.23,v_{20}=4.91\), with corresponding energy \(H=34.695\).} \label{figtraj}
\end{figure}

We can now calculate analytically the linearized Poincare map at the periodic orbits
for this hard-wall case.  By Theorem 1, this will approximate the map for the smooth impact system.  Clearly, for the smooth case, finding the periodic orbit and its stability is computationally expensive, whereas here (the hard-wall case), everything can be done analytically up to some Newton iterations.   Let the Poincare section $\Sigma$ be the positive $u_1$-axis.   By computing the eigenvalues of the return map to this section (see Appendix B for detailed calculations), we can now study the stability of the periodic orbit of the impact system for a variety of parameter values, and conclude that for sufficiently small \(\epsilon\), hyperbolic and elliptic periodic orbits will persist. 
For the impact system, we find numerically that for a fixed
\(\omega/\lambda\), the periodic orbit is hyperbolic  if \(u_{10} \in (u_{1s},u_c(\omega/\lambda,\beta,u_{1s}))\) and is elliptic if  \(u_{10} > u_c(\omega/\lambda,\beta,u_{1s})\) (we did not detect another change of stability as \(u_{10}\) is  increased further); see Figure 4 for a typical bifurcation diagram of the impact system in the \((\omega/\lambda,u_{10})\) space and the corresponding diagram for the smooth system with $\epsilon=.0001, .001, .01,.1,.2,.3$.

\begin{figure}[htbp]
\centering
  \begin{minipage}[b]{15 cm}
    \includegraphics[scale=.6]{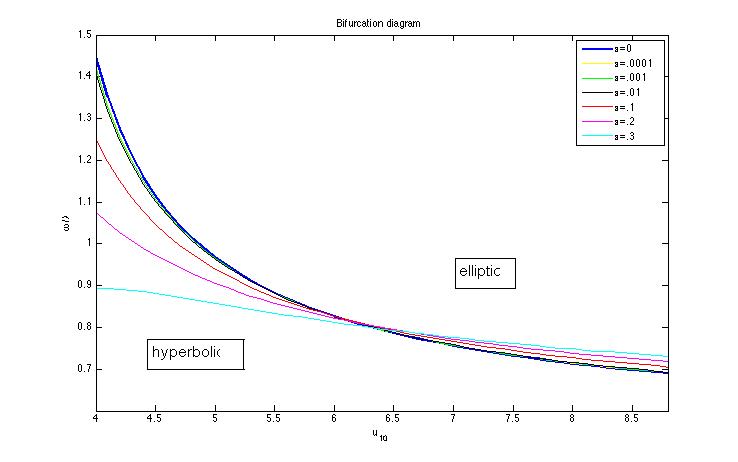}  
  \end{minipage}
\caption{Bifurcation diagram for $\epsilon=0$ (blue), $\epsilon=0.0001$ (yellow), $\epsilon=0.001$ (green),  $\epsilon=0.01$ (black), $\epsilon=0.1$ (red), $\epsilon=0.2$ (magenta), and $\epsilon=0.3$ (cyan); where $\omega=1$, $b=50$, $\beta=\pi/3$, $u_{1s}=2.5$, and $u_{2s}=0$. 
}
\label{figbif}
\end{figure}

The hyperbolic or elliptic impact trajectory may be used as an initial condition for finding the corresponding periodic orbit of the smooth impact flow by a continuation scheme in \(\epsilon\). In the continuation scheme one may choose to fix either \(u_{10}\), the energy, or the Floquet multiplier \(\Lambda\) of the periodic orbit. Figures  5 and 6 show the Poincare maps at \(\Sigma=\{u_{2}=0,v_2>0, H=H((u,v)_{periodic})\}\) for different values of $\epsilon$. In Figure 5 the Floquet multiplier  $\Lambda$ of the linearized Poincare map is kept fixed for all \(\epsilon\), and one observes that the map structure around the periodic orbit is preserved--the orbit simply shifts to the right.  In Figure 6  the energy $h$ is kept fixed for all \(\epsilon\) values.  We see that here, for some $\epsilon\in(0.1,0.2)$, the periodic orbit undergoes a bifurcation, in accordance with Figure 4.
\begin{figure}[htbp]
  \centering
  \begin{minipage}[b]{5.4 cm}
    \includegraphics[scale=.3]{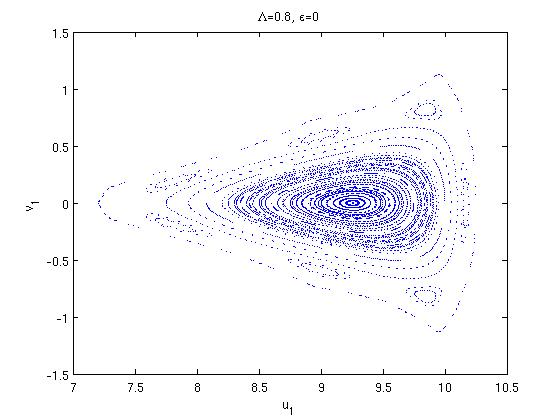}  
  \end{minipage}
  \begin{minipage}[b]{5.4 cm}
    \includegraphics[scale=.3]{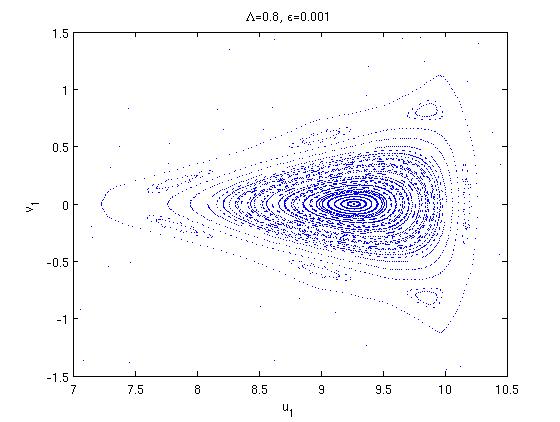}  
  \end{minipage}
  \begin{minipage}[b]{5.3 cm}
    \includegraphics[scale=.3]{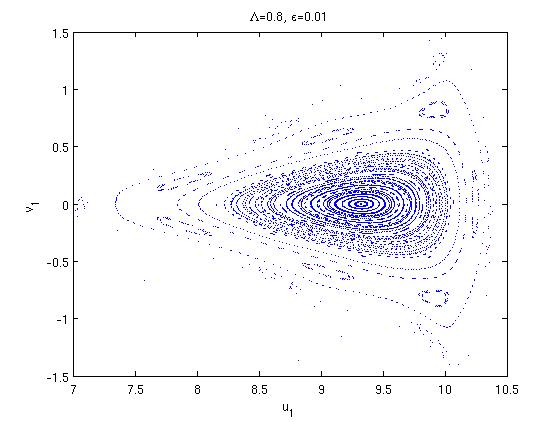}  
  \end{minipage}\\
  \begin{minipage}[b]{5.4 cm}
    \includegraphics[scale=.3]{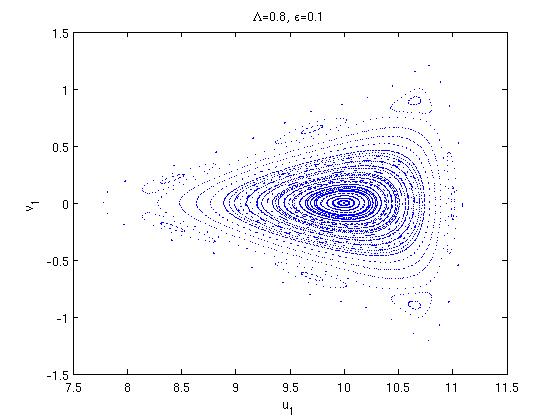}  
  \end{minipage}
  \begin{minipage}[b]{5.4 cm}
    \includegraphics[scale=.3]{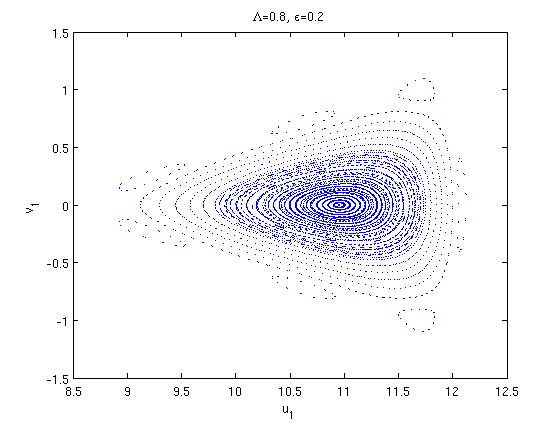}  
  \end{minipage}
  \begin{minipage}[b]{5.3 cm}
    \includegraphics[scale=.3]{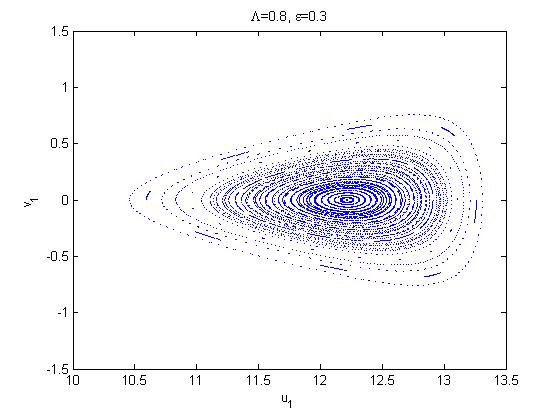}  
  \end{minipage}
  \caption{Poincare return maps for $\epsilon=0, .001, .01, .1, .2, .3$, with the periodic orbit (and thus $h$) chosen so that the real part of eigenvalues $\Lambda$ and $1/\Lambda$ of the linearized Poincare map is 0.8, using the parameters $\lambda=\sqrt{2}$, $\omega=1$, $b=50$, $\beta=\pi/3$, $u_{1s}=2.5$, $u_{2s}=0$.  In each plot, the energy is fixed according to the energy of the central periodic orbit.}
  \label{poinc_Lambdafixed}
\end{figure}

\begin{figure}[htbp]
  \centering
  \begin{minipage}[b]{5.4 cm}
    \includegraphics[scale=.27]{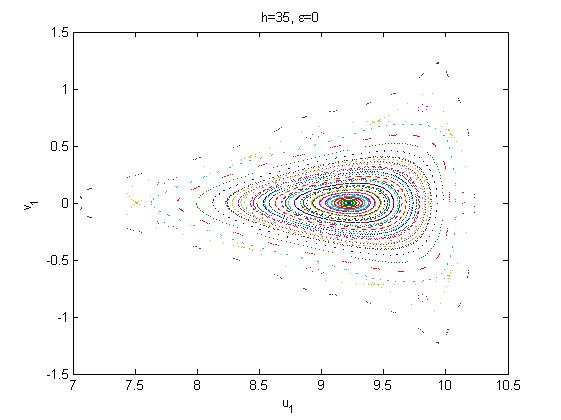}  
  \end{minipage}
  \begin{minipage}[b]{5.4 cm}
    \includegraphics[scale=.3]{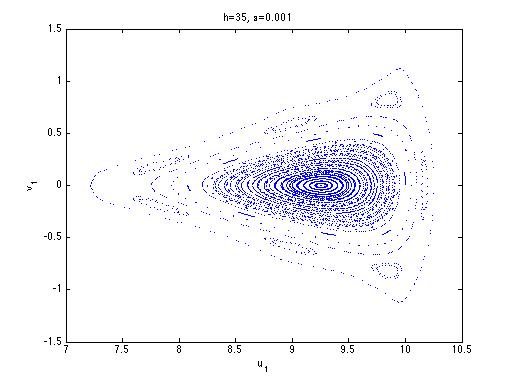}  
  \end{minipage}
  \begin{minipage}[b]{5.3 cm}
    \includegraphics[scale=.3]{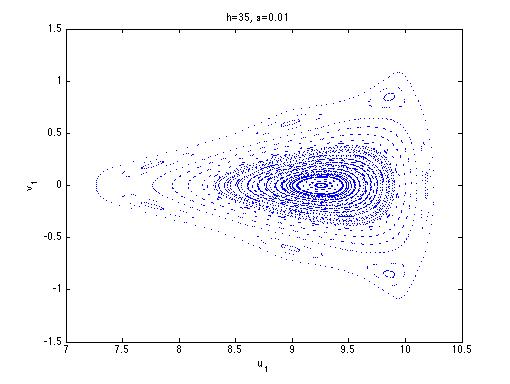}  
  \end{minipage}\\
  \begin{minipage}[b]{5.4 cm}
    \includegraphics[scale=.3]{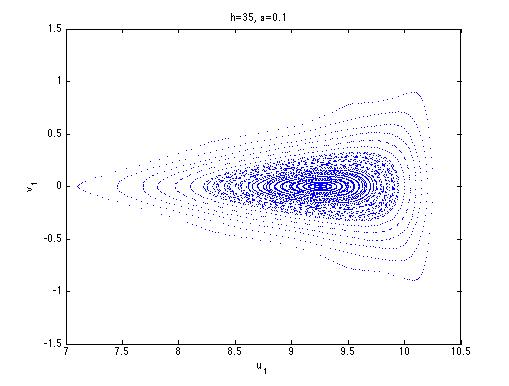}  
  \end{minipage}
  \begin{minipage}[b]{5.4 cm}
    \includegraphics[scale=.3]{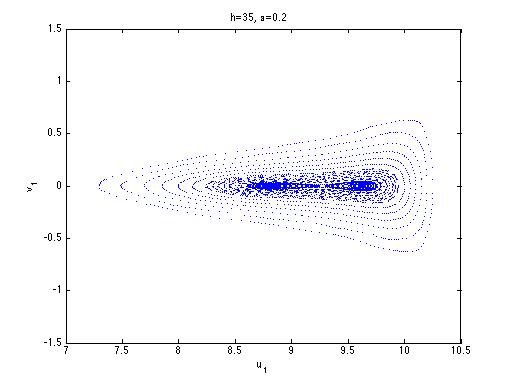}  
  \end{minipage}
  \begin{minipage}[b]{5.3 cm}
    \includegraphics[scale=.3]{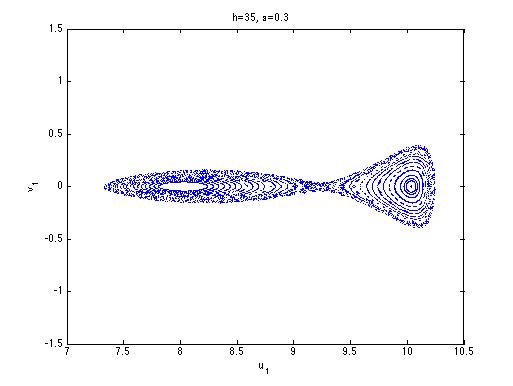}  
  \end{minipage}\\
  \caption{Poincare return maps for $\epsilon=0, .001, .01, .1, .2, .3$, with the periodic orbit chosen so that the energy is constant $h=35$, using the parameters $\lambda=\sqrt{2}$, $\omega=1$, $b=50$, $\beta=\pi/3$, $u_{1s}=2.5$, $u_{2s}=0$.  At large $\epsilon$, a change in stability is observed.}
  \label{poinc_hfixed}
\end{figure}

Finally, we demonstrate that the construction works also for the 3 d.o.f. system.  In the coupled 3 d.o.f. case  ($\delta>0$), the return map is 4d.  We project the Poincare return map to the 2d plane by restricting to the slab \begin{displaymath}
\Sigma=\{(u_1,u_2,u_3,v_1,v_2,v_3) :  u_2=0, u_3\in[-\xi, \xi], v_2>0, v_3>0, H=H   ((u,v)_{periodic})\},
\end{displaymath}
where $\xi=0.1$.  The projections of the Poincare return maps for $\delta=0$ and $\delta=0.1$ are shown in Figure 6.  Note that the blurriness in the plot for $\delta=0.1$ is due to projecting the $u_3\in[-\xi, \xi]$.  The analytical calculations for the hard-wall system were used to find an initial guess for the computations of the smooth system, leading to a dramatic decrease in computation time.  

\begin{figure}[htbp]
  \centering
  \begin{minipage}[b]{8 cm}
    \includegraphics[scale=.45]{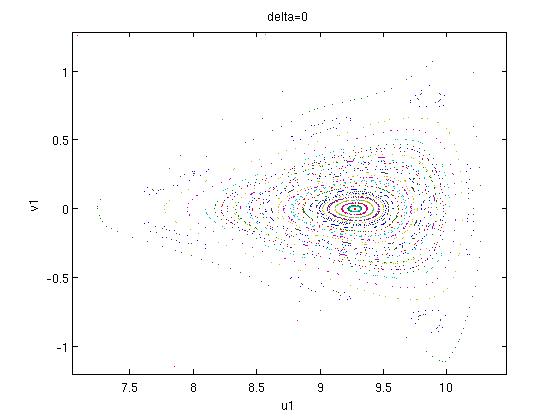}  
  \end{minipage}
  \begin{minipage}[b]{8 cm}
    \includegraphics[scale=.45]{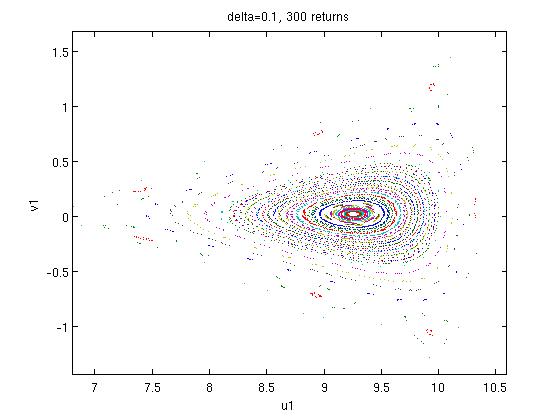}  
  \end{minipage}
  \caption{Poincare return maps for system with coupling constants $\delta=0$ and $\delta=0.1$, with periodic orbits at $u_{10}=9.27$, using the parameters $\epsilon=.001$, $\lambda=\sqrt{2}$, $\omega=1$, $b=10$, $\beta=\pi/3$, $u_{1s}=2.5$, $u_{2s}=0$, $u_{3s}=0$, $\xi=0.1$.}
  \label{poinc_coupling}
\end{figure}

Finally, we note that in this example, a linear background potential is used.  However, a similar structure would emerge if we replaced this potential with any smooth bounded separable potential having a single saddle-center fixed point (satisfying the conditions of Theorem 1)  $(U(u_1,u_2)=U_1(u_1)+U_2(u_2)$, with $U_1'(u_{1s})=0$, $U_1''(u_{1s})>0$, $U_2'(u_{2s})=0$, and $U_2''(u_{2s})<0)$.  This approximation can now be applied to a much wider class of problems, many of which are relevant to chemistry and physics, greatly reducing computation times for standard calculations such as finding periodic orbits and Poincare return maps. 

\section{Proofs of Theorems 1 and 2}
\textbf{Proof of Theorem 1}:
The proof is similar to \cite{Rapoport07}, with modifications arising when dealing with the potential term $U$.  For any given compact region $K\subset D$, the Hamiltonian flow is $C^r$-close to the impact flow  by Condition I.  Since $T$ is finite, we need only consider the behavior of the Hamiltonian flow inside a boundary layer that is close to $\partial D$.
Note that by choosing $H^*<\mathcal{E}+\hat U$ and by condition \textbf{V} we ensure that the point mass cannot escape over the barrier at \(\partial D\). Indeed, to cross the barrier   $H^*-U(q)-V(q;\epsilon)|_{\partial D}$ must be positive at some point on the boundary, but $H^*-U(q)-V(q;\epsilon)|_{\partial D}<\mathcal{E}+\hat U-\hat U-\mathcal{E}\textless0$.  

Consider an impact orbit with collision point $q_c\in\Gamma_i\setminus\Gamma^*$. Define the boundary layer near $\Gamma_i$ (where $Q(q;0)=Q_i$ by Condition IIa) as $N_\delta=\{
|Q(q;\epsilon)-Q_i|\leq \delta\}$, where $\delta$ tends to zero sufficiently slowly as $\epsilon\to 0^+$.  Take $\epsilon$ sufficiently small.  Hence, by the assumption on $\rho_T$, all collisions with the boundary occur with non-zero velocity.  The smooth impact trajectory enters $N_\delta$ at some time $t_{in}(\delta,\epsilon)$ at a point $q_{in}(\delta,\epsilon)$ close to $q_c$ with velocity $p_{in}(\delta,\epsilon)$ close to $p_0\neq 0$, and then exits $N_\delta$ at the time $t_{out}(\delta,\epsilon)$ at a point $q_{out}(\delta,\epsilon)$ with velocity $p_{out}(\delta,\epsilon)$.
Now, proving the theorem is equivalent to proving the following statements (where $k=0$ for non-degenerate tangent collisions and $k=r$ for regular collisions):
\begin{eqnarray}\label{eq8}
\lim_{\delta\to 0} \lim_{\epsilon\to 0^+}\left\| \left( q_{out}(\delta,\epsilon), t_{out}(\delta,\epsilon)\right)-\left( q_{in}(\delta,\epsilon), t_{in}(\delta,\epsilon)\right)
\right\|_{C^k}=0,
\end{eqnarray} 
which guarantees that the trajectory does not travel along the boundary, and \begin{eqnarray}\label{eq9}
\lim_{\delta\to 0} \lim_{\epsilon\to 0^+}\left\| p_{out}(\delta,\epsilon)
-p_{in}(\delta,\epsilon)+2n(q_{in})\langle p_{in}(\delta,\epsilon),n(q_{in})\rangle\right\|_{C^k}=0,
\end{eqnarray}
where $p_{out}=p_{in}-2\langle p_{in},n(q)\rangle n(q)$ and $n(q)$ is the unit inward normal to the level surface of $Q$ at the point $q$.

Without loss of generality, assume that $Q(q;0)$ increases as $q$ moves from $\partial D$ towards int$(D)$.  The partial derivatives of $Q$ satisfy \begin{eqnarray}\label{eq10}
Q_x|_{(q_c;\epsilon)}=0, \quad Q_y|_{(q_c;\epsilon)}=1.
\end{eqnarray} 
By \eqref{eq2} and Condition II, near the boundary the equations of motion have the form \begin{eqnarray}\label{eq11}
\dot x=\frac{\partial H}{\partial p_x}=p_x \quad \dot p_x=-\frac{\partial H}{\partial x}=-W'(Q;\epsilon)Q_x-U_x(x,y),
\end{eqnarray}
\begin{eqnarray}\label{eq12}
\dot y=\frac{\partial H}{\partial p_y}=p_y \quad \dot p_y=-\frac{\partial H}{\partial y}=-W'(Q;\epsilon)Q_y-U_y(x,y).
\end{eqnarray}

First, we prove \eqref{eq8} and \eqref{eq9} for $k=0$ (i.e., the $C^0$ part of the theorem, for regular or non-degenerate tangent reflections). Suppose that $\xi(\epsilon)\to 0$ sufficiently slowly, and that the orbit remains in $N_\delta$ for all $t\in I=[t_{in},t_{in}+\xi]$.  Then for all $t\in I$, 
\begin{eqnarray}\label{eq13}
q(t)=q_{in}(\delta,\epsilon)+O(\xi).
\end{eqnarray}  This follows from $p$ being uniformly bounded by $\frac{p^2}{2}=H-W(Q;\epsilon)-U\leq H^{*}-\hat U$ and from \eqref{eq11} and \eqref{eq12}.  By \eqref{eq13} and the smoothness of $U$, we have that for all $t\in I$, \begin{eqnarray}\label{eq16}
U(q(t))=U(q_{in}(\delta,\epsilon))+O(\xi).
\end{eqnarray}
And for $\xi(\epsilon)\to 0$ sufficiently slowly, \eqref{eq13} also implies that
\begin{eqnarray}\label{eq17}
q(t)=q_c+O(\xi),
\end{eqnarray}
since $q_{in}-q_c$ tends to zero as $O(\delta)$ for regular trajectories and $O(\sqrt{\delta})$ for non-degenerate tangent trajectories.  

In addition, assuming $\xi(\epsilon)\to 0$ sufficiently slowly, we claim that\footnote{This differs from \cite{Rapoport07}, in which $\alpha=1$.}, for some $\alpha\in(0,1)$,
\begin{eqnarray}\label{eq14}
p_x(t)=p_x(t_{in}(\delta,\epsilon))+O(\xi^{\alpha})
\end{eqnarray} and 
\begin{eqnarray}\label{eq15}
\frac{p_y(t)^2}{2}+W(Q(q(t);\epsilon);\epsilon)=\frac{p_y(t_{in}(\delta,\epsilon))^2}{2}
+W(\delta,\epsilon)+O(\xi^{\alpha}).
\end{eqnarray}
for all $t\in I$.  First consider \eqref{eq14}.  Note that \eqref{eq10} and \eqref{eq17} imply that 
 \begin{eqnarray}\label{eq18}
Q_x(q(t);\epsilon)=O(\xi), \quad Q_y(q(t);\epsilon)=1+O(\xi)
\end{eqnarray}
for all $t\in I$.
Given a sufficiently small $\epsilon_0$, where $\xi(\epsilon_0)\ll 1$, choose $\alpha=\alpha(\epsilon_0,U_y(q_c))\in (0,1)$ such that $\xi(\epsilon_0)^{-\alpha}=2U_y(q_c)$.  Let $K=\xi(\epsilon)^{-\alpha}$.
  It follows from the smoothness of $U$ and \eqref{eq17} that as $\epsilon\to 0$, $K\gg 2U_y(q(t))$, $t\in I$.    Divide the interval $I$\ into two sets:  $I_{<}$ where $|W'(Q;\epsilon)|<K=K(\epsilon)$ and $I_>$ where $|W'(Q;\epsilon)|\geq K=K(\epsilon)$.  In $I_<$, we have $\dot p_x=-U_x+O(\xi^{1-\alpha})$, by \eqref{eq11}, \eqref{eq13} and \eqref{eq18}.   In $I_>$, since $|W'(Q;\epsilon)|\geq K$ and $Q_y\neq 0$, we have that $\dot p_y$ is bounded away from zero, so  we can divide $\dot p_x$ in \eqref{eq11} by $\dot p_y$ in \eqref{eq12}:  \begin{displaymath}
\frac{dp_x}{dp_y}=\frac{-W'Q_x-U_x}{-W'Q_y-U_y}=\frac{Q_x+\frac{1}{W'}U_x}
{Q_y+\frac{1}{W'}U_y}=O\left(\frac{1}{K}\right)=O(\xi^\alpha).
\end{displaymath}
The change in $p_x$ on $I$ can be estimated from above as the sum of an $O(\xi)$ term (the contribution from $I_<$) plus a term of $O(\xi^{\alpha})$ times the total variation in $p_y$ (the contribution from $I_>$).  
 Recall that $p_y$ is uniformly bounded ($|p_y|\leq 2H^{*}$ from the energy constraint)
 and monotone (as $W'(Q)<-K$, $Q_y\approx 1$, and $K>\max U_y$, we have $\dot p_y>0$, see \eqref{eq12}) everywhere on $I_>$ so its total variation is indeed uniformly bounded.  This implies that the total variation in $p_y$ is uniformly bounded on $I_>$, completing the proof of \eqref{eq14}.

Now consider \eqref{eq15}.  From \eqref{eq14} and the conservation of $H=\frac{p_x^2}{2}+\frac{p_y^2}{2}+W(Q(q;\epsilon);\epsilon)+U(q)$,
we have \begin{eqnarray}\label{eq19}
\frac{p_y(t)^2}{2}+W(Q(q(t);\epsilon);\epsilon)+U(q(t))=\frac{p_y(t_{in}(\delta,\epsilon))^2}{2}
+W(\delta,\epsilon)+U(q_{in}(\delta,\epsilon))+O(\xi^{\alpha}).
\end{eqnarray}  And \eqref{eq15} follows from \eqref{eq16} and \eqref{eq19}.  

Now we claim that the time $\tau_\delta$ that the trajectory spends in the boundary layer $N_\delta$ tends to zero as $\epsilon\to 0$.  We treat the non-tangent and non-degenerate tangent cases separately.
First, consider the non-tangent trajectories (so $p_y(t_{in})$ is bounded away from zero).    
The value of $W_{in}=W_{out}=W(Q=\delta;\epsilon)$ goes to zero as $\epsilon\to 0^+$ by Condition III; thus, \eqref{eq15} implies that if $t\in I$ and  $W(Q;\epsilon)<\nu\ll \frac{p_y^2(t_{in})}{2}$, then
$p_y(t)$ is bounded away from zero.  Divide $N_\delta$ into two parts:  $N_<:=\{W:W(Q;\epsilon)\leq \nu\}$ and $N_>:=\{W:W(Q;\epsilon)>\nu\}$.    First, the trajectory enters $N_<$.  Since the value of $\displaystyle \frac{d}{dt}Q(q)=p_xQ_x+p_yQ_y$ is negative and bounded away from zero in $N_<$ (because $Q_x$ is small,  $p_y<0$ and $Q_y\approx 1$), the trajectory must reach the inner part $N_>$ by a time proportional to the width of $N_<$, which is $O(\delta)$.  And if the trajectory leaves $N_>$ after some time $t_>$, it must have $p_y>0$; thus, $t_{out}-t_{in}=O(\delta)+t_>$.  

We claim that $t_>\to 0$ as $\epsilon\to 0^+$.  We showed already that the total variation on $p_y$ is uniformly bounded.  This, along with \eqref{eq12} and Condition IV (which implies that for small $\epsilon$, $\mathcal{Q}'(W\geq \nu)\to 0$ as $\epsilon\to 0$, so $-W'(\mathcal{Q})$ is large in $N_>$ and, in particular, $-W'-U_y>0$), implies 
\begin{eqnarray*}
|t_>|\leq \frac{C}{\min_{N_>}|W'(Q;\epsilon)+U_y|}\leq\bar C\max_{N_>}|\mathcal{Q}'(W;\epsilon)|
\to 0 \quad \mbox{ as }\epsilon\to +0
\end{eqnarray*}
for some constants $C$ and $\bar C$.
Thus,  $\tau_\delta$ tends to zero as $\epsilon\to 0$ in the non-tangent case. 

Now consider the non-degenerate tangent trajectories.  First we note that the proof for the non-tangent case holds for $p_{y,in}$ tending to zero sufficiently slowly,  meaning here that $\frac{1}{Q'(W=\nu)}$ is sufficiently large with respect to $|U_y|$.  So we now prove the result for nearly-tangent trajectories, for which $p_{y,in}$ tends to zero as $\epsilon\to 0$.  First, we claim that, while the trajectory is in the boundary layer $N_\delta$ and $t-t_{in}$ is small, $p_y(t)$ remains small.  This follows from  \eqref{eq15} since 
$W(Q;\epsilon)>W_{in}=W(\delta,\epsilon)$ (since $W$ is monotone by condition IIc).  Using this and  \eqref{eq14}, we get that $p_x(t)$ is bounded away from zero as long as $q(t)\in N_\delta$ and $t-t_{in}$ is small.

Next we claim that, for a bounded away from zero interval of time starting with $t_{in}$, $\frac{d^2}{dt^2}Q(q(t);\epsilon)$ is positive and remains bounded away from zero.   $\dot Q$ is small, since  \begin{eqnarray}\label{eq20}
\dot Q:=\frac{d}{dt}Q(q(t);\epsilon)=Q_xp_x+p_yQ_y,
\end{eqnarray}
by \eqref{eq11}, \eqref{eq12}.   Recall that, by the definition of non-degenerate tangency, $p_x^TQ_{xx}p_x>U_y$ and is bounded away from zero.  We also know that $p_y$ is small, $W'(Q;\epsilon)$ is negative, $Q_y\approx 1$, and $Q_x$ is small.  Thus,  \begin{eqnarray}\label{eq21}
\frac{d^2}{dt^2}Q(q(t);\epsilon)=p_x^TQ_{xx}p_x+2Q_{xy}p_xp_y+Q_{yy}p_y^2-W'(Q;\epsilon)
(Q_x^2+Q_y^2)-U_yQ_y-U_xQ_x
\end{eqnarray}
is positive and bounded away from zero for a bounded away from zero time interval starting with $t_{in}$.

This gives, for some constant $C_1> 0$, on this time interval, 
\begin{eqnarray}\label{eq22}
Q(q(t);\epsilon)\geq Q(q_{in};\epsilon)+\dot Q(t_{in})(t-t_{in})+C_1(t-t_{in})^2.
\end{eqnarray}
We now claim that the maximum time that the nearly-tangent trajectory spends in $N_\delta=\{|Q(q;\epsilon)-Q_i|\leq \delta=|Q(q_{in};\epsilon)-Q_i|\}$ is $O(\sqrt{\delta}+p_{y,in})$ (which tends to zero).  This follows from \eqref{eq22} (which implies that this time is $O(\dot Q(t_{in}))=O(Q_x(q_{in}))+O(p_{y,in})=O(q_{in}-q_c)+O(p_{y,in})$) and from the fact that $q_{in}-q_c=O(\sqrt{\delta})$.

We have shown that, in both the non-tangent and non-degenerate tangent cases, the time the trajectory spends in the boundary layer tends to zero as $\epsilon$ goes to zero.  Thus, the proof of the theorem in the $C^0$ case is completed by 
substituting the time $\tau_\delta\to 0$ for $\xi$ in the right-hand sides of \eqref{eq13}, \eqref{eq14}, and \eqref{eq15} to get  \eqref{eq8} and \eqref{eq9} for $k=0$.

For non-tangent trajectories, we need to prove convergence in the $C^r$ topology.  Define $N_>$ and $N_<$ for small $\nu$, as in the proof of the $C^0$-convergence.  Since $\dot Q\neq 0$ in $N_<$, we can divide \eqref{eq11}, \eqref{eq12} by $\dot Q$: 
\begin{align} \label{eq23}
\frac{dq}{dQ}&=\frac{p}{Q_xp_x+p_yQ_y},\\
\label{eq30} \frac{dp}{dQ}&=\frac{-W'(Q;\epsilon)\nabla Q-\nabla U}{Q_xp_x+p_yQ_y},\\
\label{eq31} \frac{dt}{dQ}&=\frac{1}{Q_xp_x+p_yQ_y}.
\end{align}

Rewrite equations \eqref{eq23}, \eqref{eq30}, and \eqref{eq31}  in integral form:  
\begin{align}\label{eq24}
q(Q_2)-q(Q_1)&=\int_{Q_1}^{Q_2} F_q(q,p)dQ,\\
\label{eq32}
p(Q_2)-p(Q_1)&=-\int_{W(Q_1)}^{W(Q_2)} F_p(q,p)dW(Q)+\int_{Q_1}^{Q_2} F_s(q,p)dQ,\\
\label{eq33} t(Q_2)-t(Q_1)&=\int_{Q_1}^{Q_2} F_t(q,p)dQ,
\end{align}
where $F_q$, $F_p$, $F_s$, and $F_t$ denote some functions of $(q,p)$ which are uniformly bounded in the $C^r$-topology.  In (29), the first term corresponds to the integral of $\displaystyle \frac{-W'(Q)\nabla Q}{Q_x p_x+p_y Q_y}$, for which the change of variables is needed since $W'$ may be large, whereas the second term associated with the potential has regular behavior.  Note that the integrals are smal in $N_<$, since the change in $Q$ is bounded by $\delta$ and the change in $W$ is bounded by $\nu$.

Applying the successive approximation method, the Poincar\'e map (the solution to \eqref{eq24}, \eqref{eq32}, and \eqref{eq33}) from $Q=Q_1$ to $Q=Q_2$ limits to the identity map (along with all derivatives with respect to initial conditions) as $\delta,\nu\to 0$.  Thus, to prove \eqref{eq8} and \eqref{eq9} for $k=r$, we need to prove \begin{eqnarray}\label{eq25}
\lim_{\nu\to 0}\lim_{\epsilon\to 0+}\left\| ( q_{out},t_{out})-(q_{in},t_{in})
\right\|_{C^r}=0,
\end{eqnarray} \begin{eqnarray}\label{eq26}
\lim_{\nu\to 0}\lim_{\epsilon\to 0+}\left\| p_{out}-p_{in}+2n(q_{in})\langle
p_{in},n(q_{in})\rangle\right\|_{C^r}=0,
\end{eqnarray}
where $(q_{in},p_{in},t_{in})$ and $(q_{out},p_{out},t_{out})$ correspond now to the intersections of the orbit with the cross-section $W(Q(q;\epsilon),\epsilon)=\nu.$

Now we claim that $\dot p_y$ is bounded away from zero.  For any $\nu$ bounded away from zero, $\mathcal{Q}(Q;\epsilon)$ tends to zero uniformly in the $C^r$-topology as $\epsilon\to 0$ for $\nu\leq W\leq H^*$ (Condition IV).  This is also true if $\nu\to 0$ sufficiently slowly.  $W'(Q;\epsilon)=(\mathcal{Q}
'(W;\epsilon)^{-1})$ is bounded away from $K=2U_y(q_c)$ in $N_>$, so by \eqref{eq12}, the claim follows.    
Dividing \eqref{eq11} and \eqref{eq12} by \begin{displaymath}
\frac{dp_y}{dt}=-W'(Q;\epsilon)Q_y-U_y(x,y)=-(\mathcal{Q}'(W;\epsilon))^{-1}Q_y-U_y(x,y),
\end{displaymath}  we obtain
\begin{align}\label{eq27}
\frac{dq}{dp_y}&=\frac{p}{-(\mathcal{Q}'(W;\epsilon))^{-1}Q_y-U_y(x,y)}=
 \frac{\mathcal{Q}'(W;\epsilon)p}{-Q_y-\mathcal{Q}'(W;\epsilon)U_y(x,y)},\\
\label{eq34}\frac{dt}{dp_y}&=\frac{1}{-(\mathcal{Q}'(W;\epsilon))^{-1}Q_y-U_y(x,y)}=
\frac{\mathcal{Q}'(W;\epsilon)}{-Q_y-\mathcal{Q}'(W;\epsilon)U_y(x,y)},\\
\label{eq35}\frac{dp_x}{dp_y}&=\frac{-(\mathcal{Q}'(W;\epsilon))^{-1}Q_x-U_x(x,y)}
{-(\mathcal{Q}'(W;\epsilon))^{-1}Q_y-U_y(x,y)}=
\frac{-Q_x-\mathcal{Q}'(W;\epsilon)U_x(x,y)}{-Q_y-\mathcal{Q}'(W;\epsilon)U_y(x,y)},
\end{align}
where $W=H-\frac{1}{2} p^2-U(x,y)$.

By Condition IV, the $C^r$-limits as $\epsilon\to 0$ of \eqref{eq27}, \eqref{eq34}, and \eqref{eq35} are \begin{eqnarray}\label{eq28}
\frac{dq}{dp_y}=0, \quad \frac{dt}{dp_y}=0, \quad \frac{dp_x}{dp_y}=
\frac{Q_x}{Q_y}.
\end{eqnarray}
Note that \eqref{eq27}, \eqref{eq34}, and \eqref{eq35} are all bounded.  Therefore, the solution to \eqref{eq28} is the  $C^r$-limit of the solution of the system \eqref{eq27}, \eqref{eq34}, and \eqref{eq35}, since the change in $p_y$ is finite.     \eqref{eq28} implies that 
\begin{eqnarray}\label{eq40}
(q_{in},t_{in})=(q_{out},t_{out})
\end{eqnarray}
 when $\epsilon\to 0$, completing the proof to  \eqref{eq25}.  Also,    \begin{displaymath}
(p_{x,out}-p_{x,in})Q_y(q_{in};\epsilon)=(p_{y,out}-p_{y,in})Q_x(q_{in};\epsilon)
\end{displaymath}
when $\epsilon\to 0$, by \eqref{eq28}, which is equivalent to  \begin{eqnarray}\label{eq29}
p_y=\langle n(q),p\rangle, \quad p_x=p-p_yn(q).
\end{eqnarray}  \eqref{eq40} completes the proof.  $\square$\\

\begin{figure}[htbp]
\centering
  \begin{minipage}[b]{10 cm}
    \includegraphics[scale=.9]{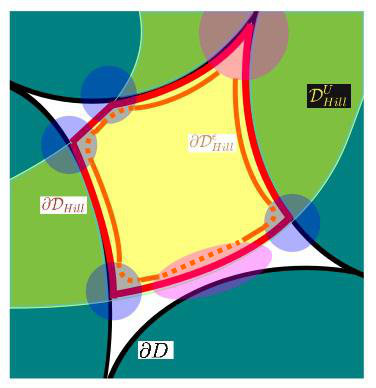}  
  \end{minipage}
  \begin{minipage}[b]{10 cm}
    \includegraphics[scale=.9]{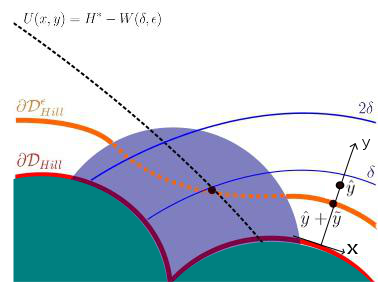}  
  \end{minipage}
\caption{Illustration of the set up for the proof of Theorem 2.  The black boundary is $\partial D$, the yellow is $\mathcal{D}^U_{Hill}(H^*)$, the red is the boundary of $\mathcal{D}_{Hill}(H^*)$, and the orange is the boundary of $\mathcal{D}_{Hill}^\epsilon(H^*)$.  The solid orange line $\partial\mathcal{D}^\epsilon_{Hill}(H^*)$ is $C^r$-close to $\partial\mathcal{D}_{Hill}(H^*)$ (in red), while the dotted orange line (inside the blue neighborhoods) is $C^0$-close.} \label{proofillustration}
\end{figure}

\noindent\textbf{Proof of Theorem 2}:\\
(i)  Recall that in the neighborhood $\tilde{N}_i$ of the billiard boundary component $\Gamma_i$  we defined  the \textit{pattern} and \textit{barrier} functions $Q_{i}(q;\epsilon),W_i(Q;\epsilon)$, and assumed (Condition II) that there exists $\epsilon_0$ such that for all $\epsilon\in (0,\epsilon_0]$ , \(V(q;\epsilon)|_{q\in\tilde{N}_i}\equiv W_i(Q(q;\epsilon)-Q_i;\epsilon)\), where \(W_{i}(Q;\epsilon)\) is monotone in this boundary layer.  In particular, there exists a \(\delta>0\) such that the thick and thin boundary layers \( N_{i}^{2\delta}(\epsilon) =\left\{ q|Q(q;\epsilon)-Q_i<2\delta,q\in \bar D\right \}\)  and \( N^{\delta}_{i}(\epsilon) =\left\{ q|Q(q;\epsilon)-Q_i<\delta,q\in \bar D\right \}\subset N^{2\delta}_{i}(\epsilon)\) are contained  inside $\tilde{N}_i$ for all $\epsilon\in (0,\epsilon_0]$; see Figure \ref{proofillustration}. Notice that the  boundary layers \(N_{i}^{2\delta}(\epsilon),N_{i}^{\delta}(\epsilon)\) are each of finite width for all \(\epsilon\) since the pattern functions have regular dependence on \(\epsilon\). Define the closed set \(K^\delta(\epsilon)=\left.\{ q|q\in D\backslash(\cup_{i}N^{\delta}_{i}(\epsilon)\cup N(\Gamma^{*}))\right.\}    \) and  let \(\rho(\epsilon)_{}=\max_{q\in K^{\delta}(\epsilon)}|V
(q;\epsilon_{})|_{C^{r+1}}\). It follows from Condition I that  \(\rho(\epsilon)\rightarrow0\) when \(\epsilon\rightarrow0\). Hence, in\footnote{This interior region does not include any of the problematic points, by definition.} $K^{\delta}(\epsilon)$,
 \(|V(q;\epsilon)|\leq \rho(\epsilon)\)  in the \(C^{r+1}\) topology and thus, in  $K^{\delta}(\epsilon)$ the level sets of \(U+V \) are \(C^{r+1}\) close to those of \(U\), as claimed. 

In each of the thick neighborhoods \(N^{2\delta}_{i}(\epsilon)\), for \(\epsilon\leq\epsilon_0\), the level sets of \(Q_{i}(q;\epsilon)\) may be viewed locally as graphs of the billiard boundary: one may set normal coordinates $(x,y)$ in each of the boundary layers, where \(x\) parameterizes the billiard boundary and  \(y\) aligns with \(\nabla Q_{}\) (hereafter we suppress the dependence on \(i\) for abbreviation). In particular, we  set \(y(\epsilon)=Q_{}(q;\epsilon)\).
 In these coordinates, \(V(x,y;\epsilon)=W(Q(x,y;\epsilon),\epsilon)=W(y(\epsilon);\epsilon)\), for all \(y\in[0,2\delta]\). Notice that \(W(y;\epsilon)\) are monotone decreasing functions, satisfying \(W(0;\epsilon)>\mathcal{E}\) (where \(W(0;\epsilon)\) may be infinite) and, by Condition\footnote{The values of \(W\) at the outer boundary may be negative, yet must converge to zero.} I,   \(|W(\delta;\epsilon)|,|W(2\delta;\epsilon)|\rightarrow0\) when \(\epsilon\rightarrow0\). 

  In   \(N^{2\delta}_{}(\epsilon)\cap  \mathcal{D}^U_{Hill}(H^*)\), away from the problematic set \(P_{Hill}(H^*)\),    \(U(x,0)<H^{*}\) (see eq. (\ref{eq:duhill})). Thus, by the monotonicity and the above observations regarding the boundary values of \(W_{}\), for such \(x\) value there exists a unique \(\hat y(\epsilon)\in(0,2\delta)\) such that \(W(\hat y(x,\epsilon);\epsilon)=H^*-U(x,0)\). Indeed, recall that in the boundary layer \(N_{i}\), for positive values of \(W\), \(W\) has an inverse \(\mathcal{Q}(W;\epsilon)\) which converges to zero along with all of its derivatives (see Condition IV). Thus, \(\hat y(x,\epsilon)=\mathcal{Q}(H^*-U(x,0);\epsilon)\) and  \(\hat y(x,\epsilon)\rightarrow 0\)  as   \(\epsilon\rightarrow0\). Expanding   \(y=\hat y(x,\epsilon)+\tilde y  \), and recalling that \(U(x,y)\) is \(C^{r+1}\) near the boundary of \(D\)   we need to find \(\tilde y(x;\epsilon)\)  such that \begin{eqnarray*}
 W(\hat y+\tilde y )=H^*-U(x,\hat y+\tilde y)=W(\hat y) -U_{y}(x,0)(\hat y+\tilde y)+O((\hat y+\tilde y)^{2}).
 \end{eqnarray*}
  Since \(\mathcal{Q'}(W(\hat y))U_{y}(x,0) (\hat y+\tilde y)\) is small, we may solve this equation by the method of successive approximation: 
\begin{displaymath}
 \tilde{y}_{j+1}=\mathcal{Q}(W(\hat{y})-U_y(x,0)(\hat{y}+\tilde{y}_j);\epsilon)-\hat{y}=
-\mathcal{Q}'(W(\hat{y});\epsilon)U_y(x,0)(\hat{y}+\tilde{y}_j)+O((\hat{y}+\tilde{y}_j)^2), \mbox{  } \tilde{y}_0=0.
\end{displaymath} 
 Thus, for any \(x\) value for which \(H^*-U(x,0)>\xi\),  the smooth and impact boundaries of the Hill's region are \(C^{r+1}\)-close, and, in particular, there exists \(\epsilon_{1}(\xi)\) such that for every such \(x\) and \(\epsilon\in(0,\epsilon_{1}]\) there exists a unique \(y(x;\epsilon)<\delta\) with  \(U(x,y(x;\epsilon))+W(y(x;\epsilon);\epsilon)=H^{*}\) and \(y(x;\epsilon)\rightarrow 0\) as \(\epsilon\rightarrow0\).  
\\

\noindent (ii)  Here we address the behavior near transverse intersections of the Hill's region corner set.   Let \((x^{*},0)\in\Sigma_{H^{*}}\ \subset\Gamma_{Hill}^*(H^*)\), namely, the level set  \(\Sigma _{U}(H^{*})=\left\{(x,y)|U(x,y)=H^*\right\}\) intersects \(\partial D\) transversely at the manifold \(\Sigma_{H^{*}}\)  parameterized by \((x^{*},0)\). It follows from the transversality assumption that \(|\nabla _{x}U(x,y)|\neq0\) in a neighborhood of \(\Sigma_{H^{*}}\). Hence, the level set  \(y=\delta\) (namely \(Q(q;\epsilon)=\delta\)) intersects \(\Sigma _{U}\) transversly along  \(\Sigma^{\delta}_{H^{*}}\)  parameterized by \((x^{*},y=\delta)\).  The transversality condition on \(U\) also guarantees that on one of the sides of   \(\Sigma^{\delta}_{H^{*}}\) the potential \(U\) becomes strictly smaller than \(H^{*}\) at some finite distance from   \(\Sigma^{\delta}_{H^{*}}\). 
 By the convergence of \(V\) to zero in \(K^{\delta}\), we know that for \(y\in[\delta,2\delta]\) the smooth Hill's region boundary is \(\rho(\epsilon)\)-close to \(\Sigma _{U}\) for these \(y \) values, so in particular, the smooth Hill's region boundary divides the surface   \(y=\delta\)   along  \(\Sigma^{\delta,\epsilon}_{H^{*}}\) (which is   \(\rho(\epsilon)\)-close to  \(\Sigma^{\delta}_{H^{*}})\)    into an interior part where \(U(x,\delta)+W(\delta;\epsilon)<H^{*}\) and an exterior part where the opposite inequality holds. This interior part, which extends to a finite neighborhood of    \(\Sigma^{\delta}_{H^{*}}\) by the above arguments, bounds the smooth Hill's region boundary near   \(\Sigma_{H^{*}}\) from the interior part of the Hill's region. The zero level set of \(Q\) (the billiard boundary, where  \(U(x,0)+W(0;\epsilon)>\hat U+\mathcal{E}>H^{*} \)), together with  the exterior part of  the surface   \(y=\delta\)    and the surface \(\Sigma _{U}(H^{*}+|W(\delta;\epsilon_{0})|)\), bounds it from its exterior side. Altogether, letting $\delta(\epsilon)$ slowly converge to zero (so that $W(\delta(\epsilon); \epsilon)\to 0$), we obtain that near a transverse corner of \(\mathcal{\partial D}_{Hill}(H^*) \), the smooth Hill's region boundary is confined to be \(C^{0}\)-close to the corner, see Figure \ref{proofillustration}.  
\\

\noindent (iii)  Now consider the case where the intersection of $\Sigma_U(H^*)$ with $\partial D$ at $\Sigma_{H^*}$ is non-transverse with $\displaystyle \left.\frac{\partial U}{\partial \bar{y}}\right|_{(x^*,0)}< 0$. Here it is more convenient to use a local Cartesian coordinate system centered at the intersection point, so that $q_b=(\bar{x},\bar{y}_b(\bar{x}))$ denotes the billiard boundary near \((x^{*},0)\). Let $\bar{y}_\eta(\bar{x};\epsilon)$ denote $\bar{y}(\bar{x;}\epsilon)$ such that $Q(\bar{x},\bar{y}(\bar{x};\epsilon),\epsilon)=\eta$. Fix $\epsilon_0>0$, $c>0$, $\xi>0$, and $\nu>0$ such that   for  all $\epsilon<\epsilon_0$,  $H^{*}<\nu +U(\bar{x}, \bar{y}_{\mathcal{Q}(\nu;\epsilon)}(\bar{x};\epsilon))<\mathcal{E}+\hat{U}$  for all $ |\bar{x}-x^*|<\xi$. Such a choice is possible: indeed, notice that   \(Q(\bar{x},\bar{y}_{\mathcal{Q}(\nu;\epsilon)}(\bar{x};\epsilon))=W^{-1}(\nu,\epsilon),  \) hence, for fixed \(\nu\) this level set of \(Q\) approaches the billiard boundary:  \(\bar{y}_{\mathcal{Q}(\nu;\epsilon)} \rightarrow\bar{y}_b(\bar{x})\). Thus, at a small neighborhood of \(q_{c}\), namely for sufficiently small  \(\xi\), as \(\epsilon\rightarrow0\), \( \ U(\bar{x}, \bar{y}_{\mathcal{Q}(\nu;\epsilon)}(\bar{x};\epsilon))=U({x^*}, 0)+O(\xi^{2},\bar{y}_{\mathcal{Q}(\nu;\epsilon)}) \) and thus for a fixed positive \(\nu\) one can always make \(\xi\) sufficiently small so the above inequality holds.

We show below that near $q_{c}$, the level set $ Q(x,y)=W^{-1}(\nu,\epsilon)$, which approaches the billiard boundary, bounds the Hill's region from the outside (on it, the energy is larger than \(H^{*}\)) whereas the level set \(Q(x,y;\epsilon)=\delta(\epsilon)\) bounds it from the inside (energy below  \(H^{*}\)). Choosing $\delta(\epsilon)\to 0$ sufficiently slowly such that $\delta(\epsilon)>W^{-1}(\nu;\epsilon)$ and $W(\delta(\epsilon);\epsilon)\to 0$ one obtains that these two level sets approach each other with the smooth Hill's region boundary in between.
  
Fix $\overline{\gamma_1}>0$ such that $\mbox{sgn}(\bar{y}_b(\bar{x}))$ is constant (or $\bar{y}_b(\bar{x})=0$) and $\bar{y}_\delta(\bar{x})>\frac{\delta}{2}$, for all $\bar{x}$ such that $|\bar{x}-x^*|<\overline{\gamma_1}$.  Choose $\gamma_2$ such that $\bar{y}_\delta(\bar{x};\epsilon)<\gamma_2$ for all $|\bar{x}-x^*|<\overline{\gamma_1}$ and all $\epsilon<\epsilon_0$.  

Let  \begin{displaymath}
\hat{\gamma_1}= 
\left\{ 
\begin{array}{cc} \sqrt{\delta}\sqrt{\frac{\left| \left.\frac{\partial U}{\partial \bar{y}}\right|_{(x^*,0)}\right|}{\left| \left.\frac{\partial^2 U}{\partial \bar{x}^2}\right|_{(x^*,0)}\right|}} & \mbox{if }\left| \left.\frac{\partial^2 U}{\partial \bar{x}^2}\right|_{(x^*,0)}\right|>c>0 \\
\xi & \mbox{otherwise}
\end{array}
\right.,
\end{displaymath}
and let $\gamma_1=\min(\overline{\gamma_1},\hat{\gamma_1},\xi)$.

For $(\bar{x},\bar{y}_{\mathcal{Q}(\nu;\epsilon)}(\bar{x};\epsilon))\in B_{\gamma_1,\gamma_2}$,  \begin{eqnarray*}
V(\bar{x},\bar{y}_{\mathcal{Q}(\nu;\epsilon)}(\bar{x};\epsilon);\epsilon)&=& V_b(\bar{x},\bar{y}_{\mathcal{Q}(\nu;\epsilon)};\epsilon)+U(\bar{x},\bar{y}_{\mathcal{Q}(\nu;\epsilon)})\\
&=& \nu +U(\bar{x},\bar{y}_{\mathcal{Q}(\nu;\epsilon)})\\
&>& H^*,
\end{eqnarray*}
by our choice of $\nu$.

For $(\bar{x},\bar{y}_\delta(\bar{x}))\in B_{\gamma_1,\gamma_2}$, recall that $\gamma_2\geq \bar{y}_\delta(\bar{x})>\frac{\delta}{2}$.   \begin{eqnarray*}
V(\bar{x},\bar{y}_\delta(\bar{x}); \epsilon)=\underbrace{V_b(\bar{x},\bar{y}_\delta(\bar{x}); \epsilon)}_{\to 0}+U(\bar{x},\bar{y}_\delta(\bar{x})).
\end{eqnarray*}
 \begin{eqnarray*}
U(\bar{x},\bar{y}_\delta(\bar{x}))&=&\underbrace{U(x^*,0)}_{=H^*}+\left.\frac{\partial U}{\partial \bar{y}}\right|_{(x^*,0)}\bar{y}_\delta(\bar{x})+\left.\frac{1}{2}\frac{\partial^2U}{\partial \bar{x}^2}\right|_{(x^*,0)}(\bar{x}-x^*)^2+h.o.t.
\end{eqnarray*}

If $\left.\frac{\partial^2 U}{\partial \bar{x}^2}\right|_{(x^*,0)}$ goes to zero, then $U(\bar{x},\bar{y}_\delta(\bar{x}))<H^*$.

Now let  $\left.\frac{\partial^2 U}{\partial \bar{x}^2}\right|_{(x^*,0)}$ be bounded away from zero.
Clearly, if $\left.\frac{\partial^2U}{\partial \bar{x}^2}\right|_{(x^*,0)}<0$, then $U(\bar{x},\bar{y}_\delta(\bar{x}))\leq H^*$.

Now suppose $\left.\frac{\partial^2U}{\partial \bar{x}^2}\right|_{(x^*,0)}>0$.  Then \begin{eqnarray*}
U(\bar{x},\bar{y}_\delta(\bar{x}))&=&H^*+\left.\frac{\partial U}{\partial \bar{y}}\right|_{(x^*,0)}\bar{y}_\delta(\bar{x})+\left.\frac{1}{2}\frac{\partial^2U}{\partial \bar{x}^2}\right|_{(x^*,0)}(\bar{x}-x^*)^2\\
&\leq& H^*+\left.\frac{\partial U}{\partial \bar{y}}\right|_{(x^*,0)}\bar{y}_\delta(\bar{x})+\left.\frac{1}{2}\frac{\partial^2U}{\partial \bar{x}^2}\right|_{(x^*,0)} \hat{\gamma_1}^2\\
&=&  H^*+\left.\frac{\partial U}{\partial \bar{y}}\right|_{(x^*,0)} \underbrace{\left(\bar{y}_\delta(\bar{x})-\frac{\delta}{2}   \right)}_{>0}\\
&<& H^*.
\end{eqnarray*}

Letting $\epsilon\to 0$ and $\delta(\epsilon)\to 0$ sufficiently slowly, we obtain that the Hill's region boundary approaches  the billiard boundary at \(q_{c}\), $C^0$-close to the singular billiard boundary.$\square$\\

 Notice that near an interior boundary point, the topology of the singular Hill's region  $\mathcal{D}_{Hill}(H)$
 does not change as \(H\) is varied through \(H^{*}\), whereas the topology does change near a bifurcating boundary point.  Thus, we expect that in the latter case a full description of the bifurcation sequence of the Hill's region boundary needs to be carried out in the \((H,\epsilon)\) plane. We leave this description to future works.
\section{Discussion}

We have extended a theorem on approximating smooth billiards with hard-wall billiards \cite{Rapoport07}, to the smooth impact case, in which free motion in the interior of the domain is replaced by motion according to a smooth bounded background potential satisfying some natural conditions.  Roughly, the theorem states that regular reflections of the smooth impact system are close to those of the hard-wall impact system.  The result allows us to, under certain conditions, approximate a smooth impact system using the limit impact system.  This is particularly useful, given the relative ease of computation for the hard-wall case.  We have applied this theorem to a geometric model for collinear triatomic chemical reactions \cite{Lerman12}, demonstrating that the simpler hard impact calculations can be used both to get qualitative information about the behavior of the smooth system, and as a tool to reduce the computation time required to find solutions in the smooth system using continuation methods.  The results presented here suggest two main future directions: singular-like results for general smooth impact systems (to be compared with those for the smooth billiard systems  \cite{Rapoport07}), and specific physical applications.

Notably, the results regarding the smooth billiard-like potentials are divided into two types: persistence-like results and singular-like results; see the recent review \cite{RomKedar12} which summarizes these works. The persistence-like results show that near regular reflections, the billiard limit and the soft steep potential are close in the \(C^{r}\) topology, hence that their local behavior near hyperbolic trajectories is similar. On the other hand, the singular-like results  show that near billiard orbits which are tangent to the boundary or go to corners, the system with the soft potential may have very different behavior than the limit system. Nonetheless, techniques for studying this singular limit by utilizing the billiard limit have been developed. Using these techniques, it was established, for example, that the soft system may have elliptic periodic orbits for arbitrarily steep potential even in cases where the hard-wall billiard is hyperbolic  \cite{Rapoport07}. Here, we explored only the persistence-like results in the soft-impact case. It turns out that this extension by itself is quite rich. Further exploration of the singular-like results for the soft impact case may provide new insights with regard to the validity and applicability of the naive impact-like system.   

In applications where the Hamiltonian is of the form (4) and the steep potential is unbounded (or $\mathcal{E} >> |U|$)  \cite{Lerman12},  one expects that at very high energies the billiard model will provide a good approximation to the dynamics \cite{Rapoport07, RomKedar12}. Here, we extend this methodology to  lower energies, where the background potential is non-negligible, yet the reflection from some boundaries is well approximated by impacts. We expect that this approach will be particularly relevant to molecular dynamics problems. There, the Pauli repulsion term is very steep and is dominant at short range at all energy levels, whereas the van der Waals or dispersion forces are smoother and contribute to a background potential that affects the motion only at energies that are of the order of the barrier energies. In \cite{Lerman12} we provided one example for this general  approach, and here, in section 3, we further explored some of the possibilities to utilize it. We believe this direction may provide qualitative insights of the dynamics in other molecular dynamics problems that are inherently non-linear and far from being integrable.  Moreover, as demonstrated, the results for the limit systems may be used to reduce computation time for calculations of various entities
such as stable periodic orbits of the smooth system, and possibly, stable and unstable manifolds. 

Finally, Theorem 1 could possibly be extended to other important cases; for example, cases in which the potential and the billiard boundary move in time  and cases in which the kinetic energy depends on the position, as when the particles are charges and are subjected to a magnetic field (see \cite{Berglund00, Berglund96, Dullin98}).  Such extension will further enhance the applicability of this methodology to additional fields of physics and chemistry.

\newpage
\appendix 
\section*{Appendix A:  Conditions for the $C^r$- and $C^0$-closeness theorems }
There are five conditions needed for proving Theorem 1. Conditions \textbf{I-IV} are concerned with the billiard-like potential \(V(q;\epsilon)\) and its limiting behavior and
are identical to those   formulated in \cite{Rapoport07,Turaev98} (repeated here in order to set up the notation for the proof). The last condition, concerned with the smooth potential \(U(q)\), is new.
   
\
   
\noindent \textbf{Condition I}:  For any fixed compact region $K\subset D$, the potential $V(q;\epsilon)$ diminishes along with all its derivatives as $\epsilon\to 0$:  \begin{displaymath}
 \lim_{\epsilon\to 0} \|V(q;\epsilon)|_{q\in K}\|_{C^{r+1}}=0.
 \end{displaymath}
 
 We assume that the level sets of $V$ may be realized by some finite function near the boundary.  Let $N(\Gamma^*)$ denote the fixed (independent of $\epsilon$) neighborhood of the corner set and $N(\Gamma_i)$ denote the fixed neighborhood of the boundary component $\Gamma_i$ in the $\mathbb{R}^d$ topology.  Define $\tilde{N}_i =N(\Gamma_i)\setminus N(\Gamma^*)$, and assume that $\tilde{N}_i\cap \tilde{N}_j=\emptyset$ when $i\neq j$.  
 
 Assume that for all small $\epsilon\geq 0$ there exists a \textit{pattern function} \begin{displaymath}
 Q(q;\epsilon):\bigcup_i \tilde{N}_i\to R^1
 \end{displaymath}
 which is $C^{r+1}$ with respect to $q$ in each of the neighborhoods $\tilde{N}_i$ and it depends continuously on $\epsilon$ (in the $C^{r+1}$ topology, so it has, along with all its derivatives, a proper limit as $\epsilon\to 0$).
 
Further assume that in each of the neighborhoods $\tilde{N}_i$ the following is fulfilled.

\noindent \textbf{Condition IIa}:  The billiard boundary is composed of level surfaces of $Q(q;0)$:  \begin{displaymath}
 Q(q;\epsilon=0)|_{q\in\Gamma_i\cap\tilde{N}_i}\equiv Q_i=\mbox{constant.}
 \end{displaymath}
 In the neighborhood $\tilde{N}_i$ of the boundary component $\Gamma_i$ (so $Q(q;\epsilon)$ is close to $Q_i$), define a \textit{barrier function} $W_i(Q;\epsilon)$, which is $C^{r+1}$-smooth in $Q$, continuous in $\epsilon$, and does not depend explicitly on $q$. Also assume that there exists $\epsilon_0$ such that the Conditions IIb-c are satisfied: 
 
 \noindent\textbf{Condition IIb}:  For all $\epsilon\in (0,\epsilon_0]$ the potential level sets in $\tilde{N}_i$ are identical to the pattern function level sets and thus \begin{displaymath}
 V(q;\epsilon)|_{q\in\tilde{N}_i}\equiv W_i(Q(q;\epsilon)-Q_i;\epsilon),
 \end{displaymath} and
 
 \noindent\textbf{Condition IIc}:  For all $\epsilon\in(0,\epsilon_0]$, $\nabla V$ does not vanish in the finite neighborhoods of the boundary surfaces $\tilde{N}_i$, thus \begin{displaymath}
 \nabla Q|_{q\in\tilde{N}_i}\neq 0
 \end{displaymath} and for all $Q(q;\epsilon)|_{q\in\tilde{N}_i}$ \begin{displaymath}
 \frac{d}{dQ} W_i(Q-Q_i;\epsilon)\neq 0.
 \end{displaymath}
 Adopt the convention that $Q>Q_i$ corresponds to the points near $\Gamma_i$ inside $D$.

\noindent \textbf{Condition III}:  There exists a  constant $\mathcal{E}_i>0$ ($\mathcal{E}_i$ may be infinite) such that as $\epsilon\to+0$ the barrier function increases from zero to $\mathcal{E}_i$ across the boundary $\Gamma_i$:
 \begin{displaymath}
  \lim_{\epsilon\to +0} W(Q;\epsilon)=\left\{ \begin{matrix}0, & Q>Q_i \\
 \mathcal{E}_i, & Q<Q_i \\
 \end{matrix}\right. .
  \end{displaymath} 
 
\noindent \textbf{Condition IV}:  As $\epsilon\to +0$, for any fixed $W_1$ and $W_2$ such that  $0<W_1<W_2<c$, for each boundary component $\Gamma_i$, the function $Q_i(W;\epsilon)$ tends to zero uniformly on the interval $[W_1,W_2]$ along with all of its $(r+1)$ derivatives.

The  last condition is concerned with the addition of the smooth component of the potential \(U(q)\), assuring that together with the billiard-like potential, particles that are initially in \(D\) cannot escape. Denoting the minimal barrier height by \(\mathcal{E}\):

\begin{equation}\mathcal{E}=\min_{i} \mathcal{E}_i\end{equation}
and the minimal value  of \(U\) on the billiard boundary
by \(\hat U\):
\begin{equation}
\hat U=\min_{q\in\partial D}U(q).
\end{equation}
To ensure that particles cannot escape from \(D\), we require that  \(\hat U>-\mathcal{E}\).  We thus impose the following condition on \(U(q)\):

\

\noindent\textbf{Condition V} :  $U(q)$ is a $C^{r+1}$-smooth potential bounded in the $C^{r+1}$ topology on an open set $\mathcal{D}$ where $\overline{D}\subset \mathcal{D}$. The minimum of \(U\) on the boundary \(\partial D\) satisfies $-\mathcal{E}<\hat U$.
 \\

\

\section*{Appendix B:  Calculation of the linearized return map}
Let the Poincare section $\Sigma$ be the positive $u_1$-axis.  Let $f$ be the flow from $\Sigma$ to the upper billiard boundary, and let $t_c$ be the collision time (which is found using a shooting method):  \begin{eqnarray*}
 f_{u_1}&:=&u_1(t_c)= u_{1s}+(u_{10}-u_{1s})\cos(\omega t_c)+\frac{v_{10}}{\omega}\sin(\omega t_c)\\
 f_{u_2}&:=& u_2(t_c)= u_{20}\cosh(\lambda t_c)+\frac{v_{20}}{\lambda}\sinh(\lambda t_c)\\
 f_{v_1}&:=& v_1(t_c)= -\omega(u_{10}-u_{1s})\sin(\omega t_c)+v_{10}\cos(\omega t_c)\\
 f_{v_2}&:=& v_2(t_c)= \lambda u_{20}\sinh(\lambda t_c)+v_{20}\cosh(\lambda t_c),
 \end{eqnarray*}
 
The choice of initial conditions fixed the energy:  \begin{eqnarray*}
H&=& \frac{v_{10}^2}{2}+\frac{v_{20}^2}{2}+\frac{\omega^2}{2}(u_{10}-u_{1s})^2-\frac{\lambda^2}{2} (u_{20}-u_{2s})^2,
\end{eqnarray*}
so \begin{eqnarray*}
v_{20}&=& \sqrt{2H-v_{10}^2-\omega^2(u_{10}-u_{1s})^2 +\lambda^2 u_{20}^2   }\\
&=&  \sqrt{2H-v_{10}^2-\omega^2(u_{10}-u_{1s})^2    }.
\end{eqnarray*}

At $t_c$ (at the upper boundary), \begin{eqnarray*}
f_{u_2}&=&f_{u_1}\tan(\beta/2)\\
f_{v_2}&=& \sqrt{2H-f_{v_1}^2-\omega^2(f_{u_1}-u_{1s})^2 +\lambda^2 f_{u_2}^2   }\\
&=&  \sqrt{2H-f_{v_1}^2-\omega^2(f_{u_1}-u_{1s})^2 +\lambda^2 f_{u_1}^2\tan(\beta/2)^2   }.
\end{eqnarray*}

Thus, we can consider instead the $2 \times 2$ system
 \begin{eqnarray*}
 f_{u_1}&:=&u_1(t_c)= u_{1s}+(u_{10}-u_{1s})\cos(\omega t_c)+\frac{v_{10}}{\omega}\sin(\omega t_c)\\
 f_{v_1}&:=& v_1(t_c)= -\omega(u_{10}-u_{1s})\sin(\omega t_c)+v_{10}\cos(\omega t_c).
\end{eqnarray*}

We calculate
\begin{displaymath}
Df|_{(u,v)_{periodic}}=D_{(u_0,v_0)}f+\frac{\partial f}{\partial t_c}\cdot \nabla_{(u_0,v_0)}t_c,
\end{displaymath}
where 
\begin{displaymath}
D_{(u_0,v_0)}f|_{(u,v)_{periodic}}=\begin{bmatrix}\cos(\omega t_c) &  \frac{1}{\omega}\sin(\omega t_c) \\
-\omega \sin(\omega t_c)  & \cos(\omega t_c)  
\end{bmatrix}
\end{displaymath}
and
\begin{displaymath}
\frac{\partial f}{\partial t_c}=\begin{bmatrix}
-\omega(u_{10}-u_{1s})\sin(\omega t_c)+v_{10}\cos(\omega t_c) \\
-\omega^2(u_{10}-u_{1s})\cos(\omega t_c)-\omega v_{10} \sin(\omega t_c)
\end{bmatrix}.
\end{displaymath}

To calculate the last matrix
\begin{eqnarray*}
\nabla_{z_0}t_c=
\begin{bmatrix}
\frac{\partial t_c}{\partial u_{10}} & \frac{\partial t_c}{\partial v_{10}}
\end{bmatrix},
\end{eqnarray*} 
we use
 \begin{eqnarray*}
 f_{u_1}:=u_1(t_c)&=& u_{1s}+(u_{10}-u_{1s})\cos(\omega t_c)+\frac{v_{10}}{\omega}\sin(\omega t_c)\\
f_{u_2}:= u_2(t_c)&=& u_{20}\cosh(\lambda t_c)+\frac{v_{20}}{\lambda}\sinh(\lambda t_c)\\
&=& \frac{\sqrt{2H-v_{10}^2-\omega^2(u_{10}-u_{1s})^2}}{\lambda} \sinh(\lambda t_c).
\end{eqnarray*}

\begin{eqnarray*}
0=F(z_0+dz_0,t_c+dt_c)&=& \tan(\frac{\beta}{2})-\frac{f_{u_2}(z_0+dz_0,t_c+dt_c)}{f_{u_1}(z_0+dz_0,t_c+dt_c)}\\
&=& \tan(\frac{\beta}{2})-\left[\frac{f_{u_2}(z_0,t_c)}{f_{u_1}(z_0,t_c)} -\frac{f_{u_2}}{f_{u_1}^2}
\left( \nabla_zf_{u_1}dz_0+\frac{\partial f_{u_1}}{\partial t_c}dt_c \right) + 
\frac{\nabla_z f_{u_2}dz_0+\frac{\partial f_{u_2}}{\partial t_c}dt_c}{f_{u_1}}  \right]\\
&=&  \underbrace{\tan(\frac{\beta}{2})-\frac{f_{u_2}(z_0,t_c)}{f_{u_1}(z_0,t_c)}}_{=0} +\frac{f_{u_2}}{f_{u_1}^2}
\left( \nabla_zf_{u_1}dz_0+\frac{\partial f_{u_1}}{\partial t_c}dt_c \right) - 
\frac{\nabla_z f_{u_2}dz_0+\frac{\partial f_{u_2}}{\partial t_c}dt_c}{f_{u_1}},
\end{eqnarray*}
so \begin{displaymath}
\frac{f_{u_2}}{f_{u_1}^2}
\left( \nabla_zf_{u_1}dz_0+\frac{\partial f_{u_1}}{\partial t_c}dt_c \right) =
\frac{\nabla_z f_{u_2}dz_0+\frac{\partial f_{u_2}}{\partial t_c}dt_c}{f_{u_1}}
\end{displaymath}
hence 
\begin{displaymath}
\tan(\frac{\beta}{2})
\left( \nabla_zf_{u_1}dz_0+\frac{\partial f_{u_1}}{\partial t_c}dt_c \right) =
\nabla_z f_{u_2}dz_0+\frac{\partial f_{u_2}}{\partial t_c}dt_c.
\end{displaymath}

Taking derivatives
\begin{eqnarray*}
 \nabla_zf_{u_1}dz_0 &=& du_{10}\cos(\omega t_c)+ dv_{10} \frac{1}{\omega}\sin(\omega t_c)\\
\frac{\partial f_{u_1}}{\partial t_c}dt_c&=& dt_c \left[ -(u_{10}-u_{1s})\omega\sin(\omega t_c) +
v_{10}\cos(\omega t_c)   \right] \\
\nabla_zf_{u_2}dz_0 &=& du_{10} \left[ \frac{-\omega^2\sinh(\lambda t_c)(u_{10}-u_{1s})}{\lambda \sqrt{2H-\omega^2(u_{10}-u_{1s})^2-v_{10}^2}}\right]  +dv_{10}\left[ \frac{-v_{10}\sinh(\lambda t_c)}{\lambda \sqrt{ 2H-v_{10}^2-\omega^2(u_{10}-u_{1s})^2} }    \right]\\
\frac{\partial f_{u_2}}{\partial t_c}dt_c &=& dt_c \cosh(\lambda t_c)\sqrt{2H-v_{10}^2-\omega^2(u_{10}-u_{1s})^2},
\end{eqnarray*}

and substituting, we get
\begin{eqnarray*}
\frac{\partial t_c}{\partial u_{10}}&=& \frac{\frac{\partial f_{u_2}}{\partial u_{10}}-\tan(\beta/2) \frac{\partial f_{u_1}}{u_{10}}  }{ \tan(\beta/2) \frac{\partial f_{u_1}}{\partial t_c}-\frac{\partial f_{u_2}}{\partial t_c}  } \\
&=& \frac{\frac{-\omega^2\sinh(\lambda t_c)(u_{10}-u_{1s})}{\lambda \sqrt{2H-\omega^2(u_{10}-u_{1s})^2-v_{10}^2}} -\tan(\beta/2)\cos(\omega t_c)
}{\tan(\beta/2)\left[ -(u_{10}-u_{1s})\omega\sin(\omega t_c) +
v_{10}\cos(\omega t_c)   \right]-\cosh(\lambda t_c)\sqrt{2H-v_{10}^2-\omega^2(u_{10}-u_{1s})^2}
}\\
\frac{\partial t_c}{\partial v_{10}}&=&  \frac{\frac{\partial f_{u_2}}{\partial v_{10}}-\tan(\beta/2) \frac{\partial f_{u_1}}{v_{10}}  }{ \tan(\beta/2) \frac{\partial f_{u_1}}{\partial t_c}-\frac{\partial f_{u_2}}{\partial t_c}  } \\
&=& \frac{ \frac{-v_{10}\sinh(\lambda t_c)}{\lambda \sqrt{ 2H-v_{10}^2-\omega^2(u_{10}-u_{1s})^2} }     -\tan(\beta/2)\frac{1}{\omega}\sin(\omega t_c)
}{\tan(\beta/2)\left[ -(u_{10}-u_{1s})\omega\sin(\omega t_c) +
v_{10}\cos(\omega t_c)   \right]-\cosh(\lambda t_c)\sqrt{2H-v_{10}^2-\omega^2(u_{10}-u_{1s})^2}}.
\end{eqnarray*}

Now let $r^{up}$ denote the reflection at the upper boundary.  By the reflection law at the upper boundary, 
 \begin{eqnarray*}
Dr^{up}=\begin{bmatrix}
1 & 0 \\
\frac{-\sin(\beta)(\omega^2(f_{u_1}-u_{1s})-f_{u_1}\lambda^2\tan(\beta/2)^2)}{
\sqrt{2H-\omega^2(f_{u_1}-u_{1s})^2-f_{v_1}^2+f_{u_1}^2\lambda^2\tan(\beta/2)^2}} &
\cos(\beta)-\frac{f_{v_1}\sin(\beta)}{\sqrt{2H-\omega^2(f_{u_1}-u_{1s})^2-f_{v_1}^2+f_{u_1}^2\lambda^2\tan(\beta/2)^2}}
\end{bmatrix}
\end{eqnarray*} with det$=1$.

Let \begin{eqnarray*}
\hat{f}_{u_1}&:=&  f_{u_1}\\
\hat{f}_{u_2}&:=&  f_{u_2}\\
\hat{f}_{v_1}&:=&  f_{v_1}\cos(\beta)+\sin(\beta) \sqrt{2H-f_{v_{1}}^2-\omega^2(f_{u_{1}}-u_{1s})^2+\lambda^2(\tan(\beta/2)f_{u_1})^2}\\
\hat{f}_{v_2}&:=&f_{v_1}\sin(\beta)-\cos(\beta) \sqrt{2H-f_{v_{1}}^2-\omega^2(f_{u_{1}}-u_{1s})^2+\lambda^2(\tan(\beta/2)f_{u_1})^2}.
\end{eqnarray*}

Let $g$ be the map from the upper boundary back to the lower boundary:  \begin{eqnarray*}
g_{u_1}&=& u_{1s}+(\hat{f}_{u_1}-u_{1s})\cos(\omega t_g)+\frac{\hat{f}_{v_1}}{\omega}\sin(\omega t_g)\\
g_{v_1}&=& -\omega (\hat{f}_{u_1}-u_{1s})\sin(\omega t_g)+\hat{f}_{v_1}\cos(\omega t_g)
\end{eqnarray*}
taking time $t_g$.  (In the periodic case, $t_g=t_c$.)

We can write \begin{eqnarray*}
g_{u_2}&=&-g_{u_1}\tan(\beta/2)\\
g_{v_2}&=& -\sqrt{2H-g_{v_1}^2-\omega^2(g_{u_1}-u_{1s})^2+\lambda^2 g_{u_2}^2}\\
&=&  -\sqrt{2H-g_{v_1}^2-\omega^2(g_{u_1}-u_{1s})^2+\lambda^2 g_{u_1}^2\tan(\beta/2)^2}.
\end{eqnarray*}

Similar to the previous calculation, 
\begin{eqnarray*}
D_{(u,v)_0}g|_{(u,v)_{periodic}}=\begin{bmatrix}
\cos(\omega  t_c)  & \frac{1}{\omega}\sin(\omega  t_c) \\
-\omega \sin(\omega t_c) & \cos(\omega  t_c)
\end{bmatrix}
\end{eqnarray*}
and
\begin{eqnarray*}
\frac{\partial g}{\partial t_g}=\begin{bmatrix}
-\omega(\hat{f}_{u_1}-u_{1s})\sin(\omega  t_c)+\hat{f}_{v_1}\cos(\omega  t_c)\\
-\omega^2(\hat{f}_{u_1}-u_{1s})\cos(\omega  t_c)-\omega \hat{f}_{v_1}\sin(\omega  t_c)
\end{bmatrix}.
\end{eqnarray*}

We calculate
\begin{displaymath}
\nabla_{z_0}t_g= \begin{bmatrix}
\frac{\partial t_g}{\partial \hat{f}_{u_1}} & \frac{\partial t_g}{\partial \hat{f}_{v_1}}
\end{bmatrix},
\end{displaymath}
using 
\begin{displaymath}
-\tan(\frac{\beta}{2})
\left( \nabla_zg_{u_1}dz_0+\frac{\partial g_{u_1}}{\partial t_c}dt_c \right) =
\nabla_z g_{u_2}dz_0+\frac{\partial g_{u_2}}{\partial t_c}dt_c.
\end{displaymath}
to get 
\begin{eqnarray*}
\frac{\partial t_c}{\partial \hat{f}_{u_1}}&=& \frac{\frac{\partial g_{u_2}}{\partial \hat{f}_{u_1}}+\tan(\beta/2)\frac{\partial g_{u_1}}{\partial \hat{f}_{u_1}}}{
-\tan(\beta/2)\frac{\partial g_{u_1}}{\partial t_c}-\frac{\partial g_{u_2}}{\partial t_c}}\\
&=& \frac{ \tan(\beta/2)\cosh(\lambda t_c)+\frac{\sinh(\lambda t_c)(\omega^2 (\hat{f}_{u_1}-u_{1s})-\hat{f}_{u_1}\lambda^2\tan(\beta/2)^2)
}{\lambda\sqrt{2H-\omega^2(\hat{f}_{u_1}-u_{1s})^2-\hat{f}_{v_1}^2+\hat{f}_{u_1}^2\lambda^2\tan(\beta/2)^2}}
  +\tan(\beta/2)  \cos(\omega t_c)
}{(**)}\\
\frac{\partial t_c}{\partial \hat{f}_{v_1}}&=& \frac{\frac{\partial g_{u_2}}{\partial \hat{f}_{v_1}}+\tan(\beta/2)\frac{\partial g_{u_1}}{\partial \hat{f}_{v_1}}}{
-\tan(\beta/2)\frac{\partial g_{u_1}}{\partial t_c}-\frac{\partial g_{u_2}}{\partial t_c}}\\
&=& \frac{\frac{\hat{f}_{v_1}\sinh(\lambda t_c)}{\lambda\sqrt{2H-\omega^2(\hat{f}_{u_1}-u_{1s})^2-\hat{f}_{v_1}^2+\hat{f}_{u_1}^2\lambda^2\tan(\beta/2)^2} }+\tan(\beta/2)\frac{1}{\omega}\sin(\omega t_c)
}{(**)}
\end{eqnarray*}

where \begin{eqnarray*}
(**)&=&-\tan(\beta/2)\left[  \hat{f}_{v_1}\cos(\omega t_c)-\omega (\hat{f}_{u_1}-u_{1s})\sin(\omega t_c)  \right] - \hat{f}_{u_1}\lambda\tan(\beta/2)\sinh(\lambda t_c)\\
& &+\cosh(\lambda t_c) \sqrt{2H-\omega^2(\hat{f}_{u_1}-u_{1s})^2-\hat{f}_{v_1}^2+\hat{f}_{u_1}^2\lambda^2\tan(\beta/2)^2}.
\end{eqnarray*}

And by the reflection law for the bottom billiard boundary, \begin{displaymath}
Dr^{low}=\begin{bmatrix}
1 & 0 \\
\frac{-\sin(\beta)\left[\omega^2(g_{u_1}-u_{1s})-g_{u_1}\lambda^2\tan(\beta/2)^2  \right]     }{\sqrt{2H-g_{v_1}^2-\omega^2( g_{u_1}-u_{1s})^2+\lambda^2(\tan(\beta/2)g_{u_1})^2} } &
\cos(\beta)-\frac{g_{v_1}\sin(\beta)}{\sqrt{2H-g_{v_1}^2-\omega^2( g_{u_1}-u_{1s})^2+\lambda^2(\tan(\beta/2)g_{u_1})^2}}
\end{bmatrix}.
\end{displaymath}

Now let \begin{eqnarray*}
\hat{g}_{u_1}&:=& g_{u_1}\\
\hat{g}_{u_2}&:=& g_{u_2}\\
\hat{g}_{v_1}&:=& g_{v_1}\cos(\beta)-g_{v_2}\sin(\beta)\\
&=&  g_{v_1}\cos(\beta)+\sin(\beta)\sqrt{2H-g_{v_1}^2-\omega^2(g_{u_1}-u_{1s})^2 +\lambda^2g_{u_2}^2}\\
&=&  g_{v_1}\cos(\beta)+\sin(\beta)\sqrt{2H-g_{v_1}^2-\omega^2(g_{u_1}-u_{1s})^2 +\lambda^2 g_{u_1}^2 \tan(\beta/2)^2}\\
\hat{g}_{v_2}&:=& -g_{v_1}\sin(\beta)-g_{v_2}\cos(\beta)\\
&=& -g_{v_1}\sin(\beta)+\cos(\beta) \sqrt{2H-g_{v_1}^2-\omega^2(g_{u_1}-u_{1s})^2 +\lambda^2 g_{u_1}^2 \tan(\beta/2)^2}
\end{eqnarray*}
and let $h$ be the map from the lower boundary back to the section, taking time $t_h$ ($t_h=t_c$ in the periodic case):
\begin{eqnarray*}
h_{u_1}&:=& u_{1s}+(\hat{g}_{u_1}-u_{1s})\cos(\omega t_c)+\frac{\hat{g}_{v_1}}{\omega}\sin(\omega t_c)\\
h_{u_2}&:=& \hat{g}_{u_2}\cosh(\lambda t_c)+\frac{\hat{g}_{v_2}}{\lambda}\sinh(\lambda t_c)\\
h_{v_1}&:=& -\omega (\hat{g}_{u_1}-u_{1s})\sin(\omega t_c)+\hat{g}_{v_1}\cos(\omega t_c)\\
h_{v_2}&:=& \lambda \hat{g}_{u_2}\sinh(\lambda t_c)+\hat{g}_{v_2}\cosh(\lambda t_c).
\end{eqnarray*}

At time $t_c$, $h_{u_2}=0$, and
\begin{eqnarray*}
h_{v_2}&=& \sqrt{2H-h_{v_1}^2-\omega^2 (h_{u_1}-u_{1s})^2},
\end{eqnarray*}
so we condisder the $2\times 2$ system \begin{eqnarray*}
h_{u_1}&:=& u_{1s}+(\hat{g}_{u_1}-u_{1s})\cos(\omega t_c)+\frac{\hat{g}_{v_1}}{\omega}\sin(\omega t_c)\\
h_{v_1}&:=& -\omega (\hat{g}_{u_1}-u_{1s})\sin(\omega t_c)+\hat{g}_{v_1}\cos(\omega t_c).
\end{eqnarray*}

The linearization of the flow back to the section (with takes time $t_c$ for a periodic orbit) is given by
\begin{displaymath}
Dh|_{(u,v)_{periodic}}=D_{(u_0,v_0)}h+\frac{\partial h}{\partial t_c}\cdot \nabla_{(u_0,v_0)} t_c,
\end{displaymath}
where
\begin{displaymath}
D_{z_0}h=\begin{bmatrix}
\cos(\omega t_c) & \frac{1}{\omega}\sin(\omega t_c)\\
-\omega \sin(\omega t_c) & \cos(\omega t_c)
\end{bmatrix},
\end{displaymath}

\begin{displaymath}
\frac{\partial h}{\partial t_c}=\begin{bmatrix}
-w(\hat{g}_{u_1}-u_{1s})\sin(\omega t_c)+\hat{g}_{v_1}\cos(\omega t_c)\\
-\omega^2 (\hat{g}_{u_1}-u_{1s})\cos(\omega t_c)-\omega\hat{g}_{v_1}\sin(\omega t_c)
\end{bmatrix},
\end{displaymath}
and 
\begin{displaymath}
\nabla_{(u_0,v_0)} t_c=\begin{bmatrix}
\frac{\partial t_c}{\partial \hat{g}_{u_1}} & \frac{\partial t_c}{\partial \hat{g}_{v_1}}
\end{bmatrix}
\end{displaymath}
with
\begin{eqnarray*}
\frac{\partial t_c}{\partial \hat{g}_{u_1}}&=& \frac{-\frac{\partial h_{u_2}}{\partial \hat{g}_{u_1}}
}{\frac{\partial h_{u_2}}{\partial t_h}}\\
&=& \frac{\tan(\beta/2)\cosh(\lambda t_c) +\frac{
\sinh(\lambda t_c)(\omega^2(\hat{g}_{u_1}-u_{1s})-\hat{g}_{u_1}\lambda^2\tan(\beta/2)^2)
}{\lambda\sqrt{2H-\hat{g}_{v_1}^2-\omega^2(\hat{g}_{u_1}-u_{1s})^2+\hat{g}_{u_1}^2\lambda^2\tan(\beta/2)^2}
}}{ \cosh(\lambda t_c)\sqrt{2H-\hat{g}_{v_1}^2-\omega^2(\hat{g}_{u_1}-u_{1s})^2+\hat{g}_{u_1}^2\lambda^2\tan(\beta/2)^2}-\hat{g}_{u_1}\lambda \tan(\beta/2)\sinh(\lambda t_c)
}\\
\frac{\partial t_c}{\partial \hat{g}_{v_1}}&=& \frac{-\frac{\partial h_{u_2}}{\partial \hat{g}_{v_1}}
}{\frac{\partial h_{u_2}}{\partial t_h}}\\
&=&\frac{\frac{\hat{g}_{v_1}\sinh(\lambda t_c)}{\lambda\sqrt{2H-\hat{g}_{v_1}^2-\omega^2(\hat{g}_{u_1}-u_{1s})^2+\hat{g}_{u_1}^2\lambda^2\tan(\beta/2)^2}
}
}{ \cosh(\lambda t_c)\sqrt{2H-\hat{g}_{v_1}^2-\omega^2(\hat{g}_{u_1}-u_{1s})^2+\hat{g}_{u_1}^2\lambda^2\tan(\beta/2)^2}-\hat{g}_{u_1}\lambda \tan(\beta/2)\sinh(\lambda t_c)
}.
\end{eqnarray*}

The linearization of the Poincare map at a periodic orbit is \begin{displaymath}
D|_{(u,v)_{periodic}}=(Dh\cdot Dr^{low}\cdot Dg \cdot Dr^{up} \cdot Df)|_{(u,v)_{periodic}}.
\end{displaymath}
 Using the Matlab symbolic math package, it was checked that it is indeed symplectic.

\
 
\noindent\textbf{Acknowledgements}\\
We acknowledge the support of the Israel Science Foundation (Grant 321/12).

\bibliographystyle{amsplain}
\bibliographystyle{plain}

\begin{thebibliography}{10}

\bibitem{Berglund00}
N.~Berglund, \emph{Classical billiards in a magnetic field and a potential},
  Nonlinear Phenomena in Complex Systems \textbf{3:1} (2000), 61--70.

\bibitem{Berglund96}
N.~Berglund and H.~Kunz, \emph{Integrability and ergodicity of classical
  billiards in a magnetic field}, Journal of Statistical Physics \textbf{83}
  (1996), 81--126.

\bibitem{BialekBook}
W.~Bialek, \emph{Biophysics: Searching for principles}, Princeton University
  Press, 2012.

\bibitem{ChampneysBook}
M.~di~Bernardo, C.J. Budd, A.R. Champneys, and P.~Kowalczyk,
  \emph{Piecewise-smooth dynamical systems: Theory and applications},
  Springer-Verlag London, 2008.

\bibitem{Dullin98}
H.R. Dullin, \emph{Linear stability in billiards with potential}, Nonlinearity
  \textbf{11} (1998), 151--173.

\bibitem{Gorelyshev08}
I.~Gorelyshev and A.~Neishtadt, \emph{On adiabadic perturbation theory for
  systems with elastic collisions}, Nonlinearity \textbf{21} (2008), 661--676.

\bibitem{KozlovBook}
V.V. Kozlov and D.V. Treschev, \emph{A genetic introduction to the dynamics of
  systems with impacts}, AMS, Providence, 1991.

\bibitem{LawleyBook}
K.P. Lawley, \emph{Advances in chemical physics, potential energy surfaces},
  John Wiley and Sons, 2009.

\bibitem{LeachBook}
A.R. Leach, \emph{Molecular modelling: Principles and applications}, Longman,
  1996.

\bibitem{Lerman12}
L.~Lerman and V.~Rom-Kedar, \emph{A saddle in a corner--a model of collinear
  triatomic reactions}, SIAM J. Appl. Dyn. Syst. (2012), to appear.

\bibitem{Lamb12}
O.~Makarenkov and J.S.W. Lamb, \emph{Dynamics and bifurcations of nonsmooth
  systems: A survey}, Physica D: Nonlinear Phenomena \textbf{241:22} (2012),
  1826–1844.

\bibitem{Rapoport07}
A.~Rapoport, V.~Rom-Kedar, and D.~Turaev, \emph{Approximating multi-dimensional
  hamiltonian flows by billiards}, Commun. Math. Phys. \textbf{272} (2007),
  567--600.

\bibitem{RomKedar12}
V.~Rom-Kedar and D.~Turaev, \emph{Billiards: a singular perturbation limit of
  smooth hamiltonian flows}, Chaos \textbf{22} (2012).

\bibitem{TannorBook}
D.J. Tannor, \emph{Introduction to quantum mechanics--a time-dependent
  perspective}, University Science Books, 2007.

\bibitem{Turaev98}
D.~Turaev and V.~Rom-Kedar, \emph{Elliptic islands appearing in near-ergodic
  flows}, Nonlinearity \textbf{11} (1998), 575--600.

\bibitem{Wu06}
T.~Wu, H.J. Werner, and U.~Manthe, \emph{Accurate pontential energy surface and
  quantum reaction rate calculations for the $\mbox{H}+\mbox{CH}_4 \to
  \mbox{H}_2+\mbox{CH}_3$ reaction}, J Chem Phys \textbf{124:16} (2006).

\end{thebibliography}

\end{document}